\documentclass[10pt]{article}
\usepackage{graphicx}
\newtheorem{theorem}{Theorem}
\newtheorem{proposition}{Proposition}
\newtheorem{lemma}{Lemma}
\newtheorem{fact}{Fact}
\newtheorem{corollary}{Corollary}

\input tcilatex

\begin{document}

\begin{titlepage}

\vskip0.5truecm

\vskip1.0truecm

\begin{center}

{\LARGE \bf Homeomorphisms of the annulus with a transitive lift}

\end{center}

\vskip  0.4truecm

\centerline {{\large F\'abio Armando Tal and Salvador Addas-Zanata}}

\vskip 0.2truecm

\centerline { {\sl Instituto de Matem\'atica e Estat\'\i stica }}

\centerline {{\sl Universidade de S\~ao Paulo}}

\centerline {{\sl Rua do Mat\~ao 1010, Cidade Universit\'aria,}} 

\centerline {{\sl 05508-090 S\~ao Paulo, SP, Brazil}}

\vskip 0.7truecm

\begin{abstract}

Let $f$ be a homeomorphism of the closed annulus $A$ that preserves 
orientation, boundary components and that has a lift $\tilde f$ to the infinite strip
$\tilde A$ which is transitive.
We show that, if the rotation number of both boundary components of $A$ is 
strictly positive, then there exists a closed nonempty connected set 
$\Gamma\subset\tilde A$ such that $\Gamma\subset]-\infty,0]\times[0,1]$, 
$\Gamma$ is unlimited,
the projection of $\Gamma$ to $A$ is dense, $\Gamma-(1,0)\subset\Gamma$ and 
$\tilde{f}(\Gamma)\subset \Gamma.$ Also, if $p_1$ is the projection in the 
first coordinate in $\tilde A$, then there exists $d>0$ such that, for any 
$\tilde z\in\Gamma,$
$$\limsup_{n\to\infty}\frac{p_1(\tilde f^n(\tilde z))-p_1(\tilde z)}{n}<-d.$$

In particular, using a result of Franks, we show that the rotation set of any 
homeomorphism of the annulus that preserves orientation, boundary components, which
 has a transitive lift without fixed points in the boundary is an interval with $0$ in its interior.
\end{abstract}

\vskip 0.3truecm

\vskip 2.0truecm


\vskip 0.8truecm

\noindent{\bf e-mail:} sazanata@ime.usp.br and fabiotal@ime.usp.br

\vskip 1.0truecm


\vfill
\hrule
\noindent{\footnotesize{Addas-Zanata and Tal are partially supported by 
CNPq, grants: 304803/2006-5 and 304360/05-8 respectively}}

\end{titlepage}

\section{Introduction and statements of the main results}

In this paper we consider homeomorphisms of the closed annulus $A=S^1\times
[0,1],$ which preserve orientation and the boundary components. Any lift of $
f$ to the universal cover of the annulus 
$\widetilde{A}=${\rm I}\negthinspace {\rm R}$\times [0,1],$ is denoted by $\widetilde{f}${\rm $,$}
a homeomorphism which satisfies $\widetilde{f}( \widetilde{x}+1,\widetilde{y}
)=\widetilde{f}(\widetilde{x},\widetilde{y})+(1,0)$ for all $(\widetilde{x}, 
\widetilde{y})\in \widetilde{A}$. 
We study properties such homeomorphisms when they have a particular lift $\widetilde{f}$ which is transitive.

Our results do not assume the existence of invariant measures of any type for
$f$, yet the importance of studying consequences of transitivity for such mappings is
underlined by the results of \cite{crovisier}, 
which imply that $\mathcal{C}^1$-generically an area preserving diffeomorphism $f$
of the closed annulus is transitive.

In order to motivate our hypotheses a little more, let us define, 
for any homeomorphism $f:A\rightarrow A$ which preserves orientation and the
boundary components and for any Borel probability $f$-invariant measure $\mu $, an invariant called the rotation number of $\mu,$
as follows:

Let $p_1:\widetilde{A}$ $\rightarrow ${\rm I}\negthinspace {\rm R} be the
projection on the first coordinate and let $p:\widetilde{A}$ $\rightarrow A$
be the covering mapping. Fixed $f$ and $\widetilde{f},$ the displacement
function $\phi :A${\rm $\rightarrow $I}\negthinspace {\rm R} is defined as

\begin{equation}
\label{defifi}\phi (x,y)=p_1\circ \widetilde{f}(\widetilde{x},\widetilde{y}%
)- \widetilde{x}, 
\end{equation}
for any $(\widetilde{x},\widetilde{y})\in p^{-1}(x,y).$ The rotation number
of $\mu $ is then given by 
$$
\rho (\mu )=\int_A\text{ }\phi (x,y)d\mu . 
$$
The importance of this definition becomes clear by Birkhoff's ergodic
theorem, which states that, for $\mu $ almost every $(x,y)$ in the annulus
and for any $(\widetilde{x},\widetilde{y})\in p^{-1}(x,y),$

$$
\rho (x,y)=\stackunder{n\rightarrow \infty }{\lim }\frac 1n\stackrel{n-1}{%
\stackunder{i=0}{\sum }}\phi \circ f^i(x,y)=\stackunder{n\rightarrow \infty 
}{\lim }\frac{p_1\circ \widetilde{f}^n(\widetilde{x},\widetilde{y})- 
\widetilde{x}}n,\text{ } 
$$
exists and 
$$
\int_A\rho (x,y)d\mu =\rho (\mu ). 
$$
Moreover, if $f$ is ergodic with respect to $\mu $, then $\rho (x,y)$ is
constant $\mu $-almost everywhere.

Following the usual definition (see \cite{alpern}), we refer to the set of
area, orientation and boundary components preserving homeomorphisms of the
annulus, which satisfy $\rho (Leb)=0$ for a certain fixed lift $\widetilde{f 
\text{ }},$ by rotationless homeomorphisms. Every time we say that $f$ is a
rotationless homeomorphism, a special lift $\widetilde{f}$ is fixed; the one
used to define $\phi .$

In \cite{alppras} it is proved that transitivity of $\widetilde{f}$ holds
for a residual subset of rotationless homeomorphisms of the annulus and the
results in 
\cite{crovisier} suggest that the same statement may hold in the $\mathcal{C}^1$ 
topology.

Our original motivation in this setting was to study a problem posed by P.
Boyland, which will be explained below.

Given a rotationless homeomorphism of the annulus $f$, by a result of Franks
(see \cite{frankscon}), if there are 2 $f$-invariant probability measures $%
\mu _1$ and $\mu _2$ with $\rho (\mu _1)<\rho (\mu _2),$ then for every
rational $\rho (\mu _1)<\frac pq<\rho (\mu _2),$ there exists a $q$-periodic
orbit for $f$ with this rotation number. So, suppose there exists a measure
with positive rotation number. By a classical result (a version of the
Conley-Zehnder theorem to the annulus) there must be fixed points with zero
rotation number, so Boyland's question is:\ Is it true that in the above
situation there must be orbits with negative rotation number? This is a very
difficult problem, which we did not solve in full generality. We considered
the following situation:

Suppose $f$ is an orientation and boundary components preserving
homeomorphism of the annulus which has a transitive lift $\widetilde{f}: 
\widetilde{A}${\rm $\rightarrow \widetilde{A}$} (one with a dense orbit). We
denote the set of such mappings by $Hom_{+}^{trans}(A).$ So every time we
say $f\in Hom_{+}^{trans}(A)$ and refer to a lift $\widetilde{f}$ of $f$, we
are always considering a transitive lift (maybe $f$ has more then one
transitive lift, we choose any of them).
Our results are the following:

\begin{theorem}
\label{teo1}:\ If $f\in Hom_{+}^{trans}(A)$ and the rotation number of $(f,
\widetilde{f})$ restricted to the boundaries of the annulus, $\rho (
\widetilde{f}\mid _{{\rm I}\negthinspace {\rm R}\times \{i\}}),$ is strictly
positive for $i=0,1,$ then there exists a closed set

$$
B^{-}\subset ]-\infty ,0]\times [0,1],B^{-}\cap \{0\}\times [0,1]\neq
\emptyset , 
$$
such that every connected component of $B^{-}$ is unlimited to the left, $
\widetilde{f}(B^{-})\subset B^{-}$ and $p(B^{-})$ is dense in the annulus.
Moreover, $B^{-}$ is the subset of 
$$
B=\stackunder{n\leq 0}{\bigcap }\widetilde{f}^n(]-\infty ,0]\times [0,1]) 
$$
which has only unlimited connected components.
\end{theorem}

Apart from the properties described in theorem 1, we also show:

\begin{theorem}
\label{teo2}: 
$$
\omega (B^{-})=\bigcap_{n=0}^\infty \overline{\left( \bigcup_{i=n}^\infty 
\widetilde{f}^i(B^{-})\right) }=\emptyset . 
$$
\end{theorem}

Thus, iterates of $B^{-}$ by $\tilde f$ converge to left end of $\widetilde{A%
}.$ The properties of $B^{-}$ allow us to extend this theorem and obtain a
stronger result:

\begin{theorem}
\label{teo3}: There exists a real number $\rho ^{+}(B^{-})<$$0$ such that,
if $\widetilde{z}\in B^{-}$, then 
$$
\stackunder{n\rightarrow \infty }{\lim \sup }\frac{p_1(\tilde f^n(\widetilde{%
z}))-p_1(\widetilde{z})}n\le \rho ^{+}(B^{-})<0. 
$$
\end{theorem}

The last theorem shows that all points in $B^{-}$ have a ``minimum negative velocity'' 
in the strip $\tilde A.$

As $\widetilde{f}$ has a dense orbit, so does $f$ and thus every point in
the annulus $A$ is non-wandering for $f.$ In this way, theorem 3 together
with Franks version of the Poincar\'e-Birkoff's theorem from \cite{annals}
implies the following:

\begin{corollary}
:\ Let $f\in Hom_{+}^{trans}(A)$. If $\tilde f$ does not have fixed points
in the boundary of $\widetilde{A}$, then the rotation set is an interval
with $0$ in its interior.
\end{corollary}

Another important consequence of theorem 3 is that, even though there
are points with rotation number in $]\rho ^{+}(B^{-}),0[,$ they do not
belong to $B^{-}.$ In particular, if such points have unstable manifolds
unbounded to the left, they must also be unbounded to the right.

Our next results, which are corollaries of the methods used to prove
theorems \ref{teo1} and \ref{teo2}, give more information on the structure
of $B^{-}.$ Their hypothesis are the same, namely, $f\in Hom_{+}^{trans}(A)$
and $\rho (\widetilde{f}\mid _{{\rm I}\negthinspace {\rm R}\times \{i\}})>0$
for $i=0,1.$

\begin{theorem}
\label{teo4}: There exists a connected component $\Gamma $ of $B^{-},$ such
that $p(\Gamma )$ is dense in the annulus, $\widetilde{f}(\Gamma )\subset
\Gamma ,$ and there is a positive integer $k$ such that $\widetilde{f}%
^{-1}(\Gamma )\subset \Gamma +(k,0),$ so $f(p(\Gamma ))=p(\Gamma ).$
\end{theorem}

\begin{theorem}
\label{teo5}: If $\Gamma $ is a connected component of $B^{-},$ then $\Gamma
-(1,0)\subset \Gamma .$
\end{theorem}

Theorems 3, 4 and 5 above have an interesting consequence. If $\Gamma $ is
as in theorem 4, and we consider the set $\Gamma _{sat}$$=\bigcup_{i=0}^%
\infty \Gamma +(i,0),$ then $\Gamma _{sat}$ is dense and connected in the
strip, $\tilde f(\Gamma _{sat})=\Gamma _{sat}$ and, in a sense, all points
in $\Gamma _{sat}$ converge to the left end of $\widetilde{A}$ through
iterations of $\tilde f$ with a strictly negative velocity. Therefore $%
\Gamma _{sat}$ can be seen as part of a dense ``unstable manifold of the
point $L$ in the $L,R$-compactification (left and right compactification) of
the strip''.

Our strategy of proof is the following: Let $f\in Hom_{+}^{trans}(A)$ be
such that $\rho (\widetilde{f}\mid _{{\rm I}\negthinspace {\rm R}\times
\{i\}})>0$ for $i=0,1.$ It is not very difficult to prove that there exists
a closed set $B^{-}\subset \left( \stackunder{n\leq 0}{\bigcap }\widetilde{f}%
^n(]-\infty ,0]\times [0,1])\right) \subset ]-\infty ,0]\times [0,1],$ such
that $B^{-}\cap \{0\}\times [0,1]\neq \emptyset ,$ every connected component
of $B^{-}$ is unlimited to the left and $\widetilde{f}(B^{-})\subset B^{-}.$
Through similar techniques, we will show both that the $\omega $-limit of $%
B^{-}$ is empty and that $p(B^{-})$ is dense in the annulus. Theorem 3 is
derived by using theorem 2 and simple properties of $B^{-}.$ The other two
theorems are proved using the machinery developed in the proof of theorem 1.

\section{Basic tools}

In this section we define several sets that will play important roles in the
proofs of the main theorems and derive some useful relations between them.

\subsection{Preliminaries}

Let us introduce the set $B^{-}$ and show some of its properties. To this
purpose, we will sometimes make use of the already mentioned left and right
compactification of $\widetilde{A}={\rm I}\negthinspace {\rm R\times [0,1],}$
denoted $L,R$-compactification, that is, we compactify the infinite strip
adding two points, $L$ (left end) and $R$ (right end), getting a closed
disk, denoted $\widehat{A}$. Clearly $\widetilde{f}$ induces a homeomorphism 
$\widehat{f}:\widehat{A}\rightarrow \widehat{A},$ such that $\widehat{f}%
(L)=L $ and $\widehat{f}(R)=R,$ see figure 1.

Given a real number $a$, let 
$$
V_a=\{a\}\times [0,1], 
$$
$$
V_a^{-}=]-\infty ,a]\times [0,1]\,{\hbox{ and}}\,V_a^{+}=[a,+\infty [\times
[0,1]. 
$$
Denote the corresponding sets on $\widehat{A}$ by $\widehat{V}_a,$ $\widehat{%
V}_a^{-}$ and $\widehat{V}_a^{+}.$ We will also denote the sets $V_0,V_0^{-}$
and $V_0^{+}$ simply by $V,V^{-}$ and $V^{+}$ respectively.

If we consider the closed set, 
$$
\widehat{B}=\stackunder{n\leq 0}{\bigcap }\widehat{f}^n(\widehat{V}^{-}), 
$$
we get that, $\widehat{f}(\widehat{B})\subset \widehat{B}$ and $L\in 
\widehat{B}.$ Denote by $\widehat{B}^{-}$ the connected component of $
\widehat{B}$ which contains $L,$ and by $B^{-}$ the corresponding set on the
strip.

\begin{lemma}
\label{bmenos}: Let $f:A\rightarrow A$ be an orientation and boundary
components preserving homeomorphism, and let $\tilde f:\widetilde{A}%
\rightarrow \widetilde{A}$ be a fixed lift of $f.$ Suppose that for every $%
a\in {\rm I}\negthinspace {\rm R,}$ there is a positive integer $n$ such
that $\tilde f^n(V)\cap V_a\not =\emptyset $ and that $\tilde f^i(V_a)\cap
V_a\not =\emptyset $ for every integer $i.$ Then $\widehat{B}^{-}\cap 
\widehat{V}\neq \emptyset $ (equivalently for the strip: $B^{-}\cap V\neq
\emptyset $)
\end{lemma}

{\it Proof:}

The proof of this result in a different context appears in Le Calvez \cite
{inst}.

Given $N>0,$ choose a sufficiently small $a<0$ such that 
$$
n=\inf \{i>0:\widetilde{f}^{-i}(V_a)\cap V\neq \emptyset \}>N. 
$$
The above is true because as $\left| a\right| $ becomes larger, it takes
more time for an iterate of $V$ to hit $V_a.$


%
%

From the definition of $n$ we get that: $\widetilde{f}^{-i}(V_a)\subset
V^{-},$ for $i=0,1,...,n-1$ and $\widetilde{f}^{-n}(V_a)\cap V\neq \emptyset
.$ This implies that there exists a simple continuous arc $\Gamma _N\subset 
\widetilde{f}^{-n}(V_a^{-})\cap V^{-},$ such that $\widehat{\Gamma }_N$
connects $L$ to $\widehat{V}$ (one endpoint of $\widehat{\Gamma }_N$ is $L$
and the other is in $\widehat{V}$), see figure 2. For this arc, if $1\leq
i\leq n,$ we get: $\widetilde{f}^i(\Gamma _N)\subset \widetilde{f}%
^{-n+i}(V_a^{-})\subset V^{-}.$ So, 
$$
\widehat{\Gamma }_N\subset \stackrel{n}{\stackunder{i=0}{\cap }}\widehat{f}%
^{-i}(\widehat{V}^{-}), 
$$
which implies, by taking the limit $N\rightarrow \infty $ $\Rightarrow $ $%
n\rightarrow \infty ,$ that $\widehat{\Gamma }_N$ has a convergent
subsequence in the Hausdorff topology to a compact connected set $\widehat{%
\Gamma }\subset \widehat{A},$ which connects $L$ to $\widehat{V}.$ From its
choice, it is clear that $\widehat{\Gamma }\subset \widehat{B}^{-}$ and thus
the lemma is proved. $\Box $

\vskip 0.2truecm

Now we:

\begin{description}
\item[Claim:]  If $\tilde f$ is transitive then the hypotheses of the
previous lemma are satisfied.
\end{description}

{\it Proof:}

The transitivity of $\widetilde{f}$ implies that we just have to prove that
for every $a\in {\rm I}\negthinspace {\rm R,}$ $\tilde f^i(V_a)\cap V_a\not
=\emptyset $ for all integers $i.$ By contradiction, suppose that for some $%
a\in {\rm I}\negthinspace {\rm R}$ and some integer $i_0,$ $\tilde
f^{i_0}(V_a)\cap V_a=\emptyset .$ Without loss of generality, we can suppose
that $\tilde f^{i_0}(V_a)\subset V_a^{-}.$ Consider the open set 
$$
W=\stackrel{i_0-1}{\stackunder{j=0}{\cup }}\tilde f^j(interior(V_a^{-})). 
$$
Clearly, $W$ is open, connected, limited to the right and $\widetilde{f}%
(W)\subset W.$ And this contradicts the existence of a dense orbit. $\Box $

\vskip 0.2truecm

Moreover, lemma \ref{bmenos} is true for any rotationless homeomorphism of
the annulus with rotation interval not reduced to zero.

\begin{proposition}
:\ If $(f,\widetilde{f})$ is a rotationless homeomorphism such that given $%
M>0,$ there exists an integer $n>0$ and a point $\widetilde{z}\in
[0,1]\times [0,1]$ such that $\left| p_1(\widetilde{f}^n(\widetilde{z}%
))\right| >M,$ then $\widehat{B}^{-}\cap \widehat{V}\neq \emptyset .$
\end{proposition}

{\it Proof:}

From lemma \ref{bmenos} it suffices to prove that for every real number $a$
there is a positive $n$ such that $\tilde f^n(V)\cap V_a\not =\emptyset ,$
because since $(f,\widetilde{f})$ is rotationless, for all real $a,$ $\tilde
f^i(V_a)\cap V_a\neq \emptyset ,$ for all integers $i.$

Suppose by contradiction that for some real $b,$ $\widetilde{f}%
^{-i}(V_b)\cap V=\emptyset $ for all integers $i>0.$ As we said above, $\rho
(Leb)=0$ implies that $\widetilde{f}^l(V_b)\cap V_b\neq \emptyset $ for all
integers $l$, so, if we suppose that $b<0,$ then $\widetilde{f}%
^{-i}(V_b)\subset V^{-}$ for all $i>0.$ And this implies that 
$$
\stackunder{n\geq 0}{\cup }\widetilde{f}^{-n}(int(V_b^{-}))\subset V^{-}. 
$$
As $\widetilde{f}^{-1}\left( \stackunder{n\geq 0}{\cup }\widetilde{f}%
^{-n}(int(V_b^{-}))\right) \subset \stackunder{n\geq 0}{\cup }\widetilde{f}%
^{-n}(int(V_b^{-})),$ there is a boundary component of the open connected
set $\stackunder{n\geq 0}{\cup }\widetilde{f}^{-n}(int(V_b^{-})),$ denoted $%
K,$ which is compact connected and intersects ${\rm I}\negthinspace {\rm %
R\times \{0\}}$ and ${\rm I}\negthinspace {\rm R\times \{1\}.}$ Clearly, $
\widetilde{f}^{-1}(K)\subset K\cup \left( \stackunder{n\geq 0}{\cup } 
\widetilde{f}^{-n}(V_b^{-})\right) .$ From $\rho (Leb)=0,$ we get that 
$$
\widetilde{f}^{-1}(K)\subset K\Rightarrow \widetilde{f}^{-1}(K+(s,0))\subset
K+(s,0) 
$$
for all integers $s$, something that contradicts the proposition hypotheses.
If $b>0,$ an analogous argument using $\widetilde{f}^{-i}(V_b)\subset V^{+}$
for all $i>0$ works. $\Box $

\vskip 0.2truecm

In the rest of the paper we assume that $f\in Hom_{+}^{trans}(A)$ and $\rho
( \widetilde{f}\mid _{{\rm I}\negthinspace {\rm R}\times \{i\}})>0$ for $%
i=0,1. $ So, from lemma \ref{bmenos}, we know that $B^{-}\subset \widetilde{A%
}$ is a closed set, limited to the right ($B^{-}\subset V^{-}$), whose
connected components (which may be unique) are all unlimited to the left,
and at least one connected component of $B^{-}$ intersects $V.$

An important point here is that, as the rotation numbers in the boundaries
of the annulus are both positive, $B$ and thus $B^{-},$ do not intersect $%
{\rm I}\negthinspace 
{\rm R\times \{0\}}$ and ${\rm I}\negthinspace {\rm R\times \{1\}}$ (because 
$\widetilde{f}(B)\subset B\subset V^{-}$). So the only part of theorem \ref
{teo1} that still has to be proved is: $p(B^{-})$ is dense in the annulus.

\subsection{The limit set of $B^{-}$}

In this subsection we examine some properties of the set 
\begin{equation}
\label{defomelim}\omega (\widehat{B}^{-})=\bigcap_{n=0}^\infty \overline{%
\bigcup_{i=n}^\infty \widehat{f}^i(\widehat{B}^{-})}, 
\end{equation}
a subset of $\widehat{A},$ and the corresponding set $\omega (B^{-})\subset 
\widetilde{A}.$

Since $\widehat{f}(\widehat{B}^{-})\subset \widehat{B}^{-},$ and since $
\widehat{B}^{-}$ is closed, we have 
$$
\omega (\widehat{B}^{-})=\bigcap_{n=0}^\infty \widehat{f}^n(\widehat{B}%
^{-}), 
$$
therefore $\omega (\widehat{B}^{-})$ is the intersection of a nested
sequence of compact connected sets, and so it is also a compact connected
set. Moreover, definition (\ref{defomelim}) implies the following lemma:

\begin{lemma}
\label{prop_omega_lim}: $\omega (B^{-})$ is a closed, $\widetilde{f}$%
-invariant set, whose connected components are all unbounded.
\end{lemma}

{\it Proof:}

Since $L\in \widehat{B}^{-}$ and $\widehat{f}(L)=L,$ we get that $L\in
\omega (\widehat{B}^{-}).$ This implies, since $\omega (\widehat{B}^{-})$ is
connected, that each connected component of $\omega (B^{-})$ is unbounded.
The other properties follow directly from the previous considerations. $\Box 
$

\vskip 0.2truecm

Of course, since $B^{-}$ is closed, we also have that $\omega (B^{-})\subset
B^{-},$ and as such, $\omega (B^{-})\cap {\rm I}\negthinspace {\rm R}\times
\{i\}=\emptyset ,\,i\in \{0,1\},$ and $\omega (B^{-})\subset V^{-}$. It is
still possible that $\omega (B^{-})=\emptyset ,$ and this is in fact true as
we show later, but for the moment we can make use of the fact that both $%
B^{-}$ and $\omega (B^{-})$ have similar properties to shorten our proofs.
For this, let $D\subset \widetilde{A}$ be a non-empty closed set with the
following properties:

\begin{itemize}
\item  {$\widetilde{f}(D)\subset D;$}

\item  {$D\subset V^{-};$}

\item  {Every connected component of $D$ is unbounded;}

\item  {$D\cap {\rm I}\negthinspace {\rm R}\times \{i\}=\emptyset ,\,i\in
\{0,1\};$}

\item  {If $\widetilde{z}\in D$ then $\widetilde{z}-(1,0)\in D.$ }
\end{itemize}

It is easily verified that $B^{-}$ has these properties, as does $\omega
(B^{-})$ if it is nonempty, so every result shown for $D$ must hold in the
particular cases of interest for us. Later, in the proof of theorem 4, we
find another set with the properties listed above.

\subsection{On the structure of $p(D)\subset A$}

First, let us start with the following lemma:

\begin{lemma}
\label{conts1}: $\overline{p(D)}\supset S^1\times \{0\},$ or $\overline{p(D)}%
\supset S^1\times \{1\}.$
\end{lemma}

{\it Proof: }

Suppose that lemma is false. Then, there are points $P_0\in S^1\times \{0\}$
and $P_1\in S^1\times \{1\}$ such that $\{P_0,P_1\}\cap $ $\overline{p(D)}%
=\emptyset .$ As $\left( \overline{p(D)}\right) ^c$ is an open set, there
exists $\epsilon >0$ such that $B_\epsilon (P_i)\cap \overline{p(D)}%
=\emptyset ,$ for $i=0,1.$ As 
$$
\widetilde{f}(D)\subset D\text{ and }p\circ \widetilde{f}(\widetilde{x}, 
\widetilde{y})=f\circ p(\widetilde{x},\widetilde{y}) 
$$
we get that 
$$
f(p(D))\subset p(D)\text{ }\Rightarrow \text{ }f(\overline{p(D)})\subset 
\overline{p(D)}. 
$$

Since $\widetilde f$ is transitive, it follows that $f$ is also transitive
and so there exists $N>0$ such that $f^{-N}(B_\epsilon (P_0))\cap B_\epsilon
(P_1)\neq \emptyset .$

Now, we must have that 
$$
f^{-N}(B_\epsilon (P_0))\cap \overline{p(D)}=\emptyset , 
$$
for if this was not true, it would imply $B_\epsilon (P_0)\cap f^N(\overline{%
p(D)})\neq \emptyset ,$ which, in turn, would imply $B_\epsilon (P_0)\cap 
\overline{p(D)}\neq \emptyset ,$ because $f^N(\overline{p(D)})\subset 
\overline{p(D)},$ contradicting the choice of $\epsilon >0.$

\vskip 0.2truecm

As $f^{-N}(B_\epsilon (P_0))\cup B_\epsilon (P_1)$ is disjoint from $
\overline{p(D)},$ this implies that there exists a simple continuous arc $%
\gamma \subset \left( f^{-N}(B_\epsilon (P_0))\cup B_\epsilon (P_1)\right) ,$
disjoint from $\overline{p(D)},$ such that its endpoints are in $S^1\times
\{0\}$ and in $S^1\times \{1\}$ $.$ So, if $\widetilde{\gamma }$ is a
connected component of $p^{-1}(\gamma ),$ then $\widetilde{\gamma }%
-(i,0)\cap D=\emptyset $ for all integers $i>0.$ 
And this contradicts the fact that the connected components of $D$ 
are unlimited to the left. $\Box $

\vskip 0.2truecm

\begin{lemma}
\label{fullmeas} If $\overline{p(D)}\neq A,$ then $\overline{p(D)}^c$ has a
single connected component which is dense in $A.$ Moreover, $\overline{p(D)}%
^c$ contains a homotopically non trivial simple closed curve in the open
annulus $S^1\times ]0,1[.$
\end{lemma}

{\it Proof:}

As $f(\overline{p(D)})\subset \overline{p(D)},$ we get that $f(\left( 
\overline{p(D)}\right) ^c)\supset \left( \overline{p(D)}\right) ^c,$ which
implies that $f^{-1}(\left( \overline{p(D)}\right) ^c)\subset \left( 
\overline{p(D)}\right) ^c.$ %
%
If $\left( \overline{p(D)}\right) ^c$ is not dense, then as $f$ is
transitive, there exists an open ball $U\subset A,$ $U\cap \left( \overline{%
p(D)}\right) ^c=\emptyset $ and a point $z\in U$ and an integer $n>0$ such
that $f^n(z)\in \left( \overline{p(D)}\right) ^c.$ And this contradicts the
fact that $f^{-1}(\left( \overline{p(D)}\right) ^c)\subset \left( \overline{%
p(D)}\right) ^c,$ so $\left( \overline{p(D)}\right) ^c$ is dense.

Let $E$ be a connected component of $\left( \overline{p(D)}\right) ^c.$
Assume by contradiction that there is no simple closed curve $\gamma \subset
E$ which is homotopically non trivial as a curve of the annulus.

In this case, $p^{-1}(E)$ is not connected and there exists an open
connected set $E_{lift}\subset \widetilde{A},$ such that $E_{lift}\cap
E_{lift}+(i,0)=\emptyset ,$ for all integers $i\neq 0$ and 
$$
p^{-1}(E)=\stackrel{+\infty }{\stackunder{i=-\infty }{\cup }}\left(
E_{lift}+(i,0)\right) . 
$$

As $f$ is transitive and $f^{-1}(\left( \overline{p(D)}\right) ^c)\subset
\left( \overline{p(D)}\right) ^c$, there exists a first $N>0,$ such that $%
f^{-N}(E)\subset E$ (in other words, for $i\in \{1,2,...,N-1\},$ $%
f^{-i}(E)\cap E=\emptyset $)$.$ This means that 
\begin{equation}
\label{usaorbdensa} 
\begin{array}{c}
\widetilde{f}^{-i}(E_{lift})\cap p^{-1}(E)=\emptyset \text{ for all }i>0 
\text{ which is not a} \\ \text{ multiple of }N\text{ and }\widetilde{f}%
^{-N}(E_{lift})\subset E_{lift}+(i_0,0), \\ \text{ for some fixed integer }%
i_0,\text{ which implies that} \\ \text{ }\widetilde{f}^{-k.N}(E_{lift})%
\subset E_{lift}+(k.i_0,0),\text{ for all integers }k>0. 
\end{array}
\end{equation}

Suppose $i_0\geq 0.$ As $\widetilde{f}$ has a dense orbit, there exists a
point $\widetilde{z}\in E_{lift}-(1,0)$ such that $\widetilde{f}^l( 
\widetilde{z})\in E_{lift},$ for some $l>0.$ But this means that, $
\widetilde{f}^{-l}(E_{lift})\cap E_{lift}-(1,0)\neq \emptyset ,$ something
that contradicts (\ref{usaorbdensa}), because we assumed that $i_0\geq 0.$ A
similar argument implies that $i_0$ can not be smaller then zero. Therefore,
all connected components of $\left( \overline{p(D)}\right) ^c$ contain a
homotopically non trivial simple closed curve.

Let $E$ be a connected component of $\left( \overline{p(D)}\right) ^c$ and $%
\gamma _E\subset E$ be a homotopically non trivial simple closed curve.
Since $f$ is transitive, $f^{-1}(\gamma _E)\cap \gamma _E\neq \emptyset .$
Thus $f^{-1}(E)\cap E\neq \emptyset $ and so, since $f^{-1}(\left( \overline{%
p(D)}\right) ^c)\subset \left( \overline{p(D)}\right) ^c$, $f^{-1}(E)\subset
E.$ But $E$ is open and $f$ is transitive, so $E$ is dense and therefore it
is the only connected component of $\left( \overline{p(D)}\right) ^c,$
proving the lemma. $\Box $

\vskip 0.2truecm

So, let 
\begin{equation}
\label{defgamaE}\gamma _E\subset \left(\overline{p(D)}\right)^c\cap
interior(A) 
\end{equation}
be a homotopically non trivial simple closed curve and let $\gamma
_E^{-}\supset S^1\times \{0\}$ and $\gamma _E^{+}\supset S^1\times \{1\}$ be
the open connected components of $\gamma _E^c.$ As $\overline{p(D)}\cap
\gamma _E=\emptyset ,$ we obtain that $\overline{p(D)}\subset \gamma
_E^{-}\cup \gamma _E^{+}.$ %
%
%
%
%
%
%
%
%
%
%
%

\begin{lemma}
\label{maisimp}:\ Let $\Gamma $ be a connected component of $D.$ If $
\overline{p(D)}\not =A$, we have:

$$
\left\{ 
\begin{array}{c}
\text{if }\overline{p(\Gamma )}\subset \gamma _E^{-},\text{ then }\overline{%
p(\Gamma )}\supset S^1\times \{0\} \\ \text{if }\overline{p(\Gamma )}\subset
\gamma _E^{+},\text{ then }\overline{p(\Gamma )}\supset S^1\times \{1\} 
\end{array}
\right. 
$$
\end{lemma}

{\it Proof: }

First note that $\overline{p(\Gamma )}\subset \overline{p(D)},$ which
implies that $\overline{p(\Gamma )}^c\supset \overline{p(D)}^c=E.$ As $E$ is
open, connected and dense, and since $\overline{p(\Gamma )}^c$ is also open,
every connected component of $\overline{p(\Gamma )}^c$ contains a point of $%
E.$ Therefore $\overline{p(\Gamma )}^c$ is also an open connected dense
subset of the annulus. Without loss of generality, suppose that $\overline{%
p(\Gamma )}\subset \gamma _E^{-}.$ This implies that $\overline{p(\Gamma )}%
\cap S^1\times \{1\}=\emptyset .$ If $\overline{p(\Gamma )}$ does not
contain $S^1\times \{0\},$ then there exists a simple continuous arc $%
\lambda $ in the annulus, which avoids $\overline{p(\Gamma )}$ and connects
some point $P_0\in \left( S^1\times \{0\}\right) \backslash \overline{%
p(\Gamma )}$ to some point $P_1\in S^1\times \{1\}.$ This is true because $%
P_0,P_1\in \overline{p(\Gamma )}^c,$ which is an open connected set. But
this means that 
$$
p^{-1}(\lambda )\cap \Gamma =\emptyset , 
$$
and this is a contradiction because each connected component of $%
p^{-1}(\lambda )$ is compact and $\Gamma $ is connected and unlimited to the
left. So, $\overline{p(\Gamma )}\supset S^1\times \{0\}.$ The other
possibility ($\overline{p(\Gamma )}\subset \gamma _E^{+}$) is held in a
similar way. $\Box $

\vskip 0.2truecm

Without loss of generality we can suppose that there exists a connected
component $\Gamma $ of $D$ that satisfies: $p(\Gamma )\subset \gamma _E^{-}.$
Thus, lemma \ref{maisimp} implies the following fact, which is the most
important information of this subsection:

\begin{fact}
\label{importante}: If $\overline{p(D)}\neq A,$ then there exists a
connected component $\Gamma $ of $D$ that satisfies: $p(\Gamma )\subset
\gamma _E^{-},$ $\overline{p(\Gamma )}\supset S^1\times \{0\},$ and thus,
given $\epsilon >0,$ $p(\Gamma )\cap S^1\times [0,\epsilon ]\neq \emptyset .$
\end{fact}

{\it Proof: }

Immediate. $\Box $

\vskip 0.2truecm

\subsection{On the structure of $D\subset \widetilde{A}$}

Let $\Gamma $ be a connected component of $D.$ We recall that, by the
definition of $D,$ $\Gamma $ is unlimited to the left.

The next proposition is used in several arguments in the remainder of the
paper.

\begin{proposition}
\label{compconex}: $\Gamma ^c$ has only one connected component, which is
unlimited.
\end{proposition}

{\it Proof:}

Clearly, there is one connected component of $\Gamma ^c$ which contains $%
int(V^{+}),$ ${\rm I}\negthinspace {\rm R\times \{0\}}$ and ${\rm I}%
\negthinspace {\rm R\times \{1\}.}$

So, if by contradiction, we suppose that $\Gamma ^c$ has another connected
component, denoted $C$, contained in $V^{-},$ its boundary must be contained
in $\Gamma .$ As $\widetilde{f}^n(\Gamma )\subset V^{-}$ for all $n\geq 0,$
we get that $\widetilde{f}^n(C)\subset V^{-}$ for all $n\geq 0.$ So for
every $\widetilde{z}\in C,$ $\stackunder{n\rightarrow \infty }{\lim \sup }$ $%
p_1\circ \widetilde{f}^n(\widetilde{z})<0.$ As $C$ is an open subset of $
\widetilde{A},$ this contradicts the transitivity of $\widetilde{f}$. $\Box $

\vskip 0.2truecm

For a connected component $\Gamma $ of $D,$ let us define 
\begin{equation}
\label{mgama}m_\Gamma =\sup \{\widetilde{x}\in {\rm I}\negthinspace {\rm R:}%
( \widetilde{x},\widetilde{y})\in \Gamma ,\text{ for some }\widetilde{y}\in
[0,1]\}\leq 0. 
\end{equation}

Consider the connected closed set $\Gamma \cup \{m_\Gamma \}\times [0,1].$
Its complement has two open unbounded connected components in%
$$
]-\infty ,m_\Gamma [\times [0,1], 
$$
one of which contains $]-\infty ,m_\Gamma [\times \{0\}$ (denoted $\Gamma
_{down}$) and another one which contains $]-\infty ,m_\Gamma [\times \{1\}$
(denoted $\Gamma _{up}$). It is possible that $\left( \Gamma \cup \{m_\Gamma
\}\times [0,1]\right) ^c$ has other unbounded connected components. But only 
$\Gamma _{up}$ and $\Gamma _{down}$ will be of interest to us, because of
the following fact, whose proof is an exercise which depends only on the
connectivity of $\Gamma $ (see lemma \ref{partecha2} for a generalization of
this result): \ 

\begin{fact}
\label{comebola}: Given a connected component $\Gamma $ of $D,$ if $\Theta $
is a closed connected set, unlimited to the left, which satisfies $\Theta
\cap \Gamma =\emptyset $ and $\Theta \subset $ $]-\infty ,m_\Gamma [\times
[0,1],$ then $\Theta \subset \Gamma _{up}$ or $\Theta \subset \Gamma
_{down}. $
\end{fact}

\vskip 0.2truecm

In the following, we will generalize the above construction and present some
simple results on the connected components of $D.$ These results will permit
us to define an order $\prec $ on the connected components of $D.$ Moreover,
it will be clear that any two disjoint closed unlimited connected sets $%
\Theta _1,\Theta _2\subset V^{-},$ which have connected complements will be
related by this order, that is either $\Theta _1\prec \Theta _2$ or $\Theta
_2\prec \Theta _1.$ This will be of importance to us, because, if $\Gamma
_1,\Gamma _2$ are connected components of $D,$ then $\widetilde{f}(\Gamma
_1) $ and $\widetilde{f}(\Gamma _2)$ may not be, they are just contained in
connected components of $D.$ As will be explained below, it is possible that
a single connected component of $D,$ denoted $\Theta ,$ contains $\widetilde{%
f}(\Gamma _1)$ and $\widetilde{f}(\Gamma _2)$ even when $\Gamma _1\cap
\Gamma _2=\emptyset $ (remember that as $\Gamma _1,\Gamma _2$ are connected
components of $D,$ either $\Gamma _1=\Gamma _2$ or $\Gamma _1\cap \Gamma
_2=\emptyset $). But in this case, if $\Gamma _1\prec \Gamma _2$ ($\Gamma
_2\prec \Gamma _1$), as $\widetilde{f}(\Gamma _1)$ and $\widetilde{f}(\Gamma
_2)$ are disjoint closed unlimited connected sets which have connected
complements, we will show that $\widetilde{f}(\Gamma _1)\prec \widetilde{f}%
(\Gamma _2)$ ($\widetilde{f}(\Gamma _2)\prec \widetilde{f}(\Gamma _1)$),
that is $\widetilde{f}$ preserves the order.

To begin, let $\Gamma \subset D$ be a connected component and let $a\leq 0$
be such that $V_a$ intersects $\Gamma .$

Consider the following open set, 
\begin{equation}
\label{defgammacompa} \Gamma ^{comp,a}=\Gamma^c \cap ]-\infty ,a[\times
[0,1]. 
\end{equation}

\begin{lemma}
\label{partecha1}: $\Gamma ^{comp,a}$ has at least two (open) connected
components, one denoted $\Gamma _{a,down}^{\prime }$ which contains $%
]-\infty ,a[\times \{0\}$ and one denoted $\Gamma _{a,up}^{\prime }$ which
contains $]-\infty ,a[\times \{1\}.$
\end{lemma}

{\it Proof:}

Suppose the lemma is false. Then, there exists $P\in ]-\infty ,a[\times
\{0\},$ $Q\in ]-\infty ,a[\times \{1\}$ and a simple continuous arc $\eta
\subset \Gamma ^{comp,a},$ whose endpoints are $P,Q.$ Clearly, $\eta \subset
]-\infty ,a[\times [0,1].$ As $\eta \cap \Gamma =\emptyset ,$ $\Gamma $ is
unlimited to the left and $\Gamma $ intersects $V_a,$ we obtain that $\Gamma 
$ intersects both connected components of $\eta ^c,$ something that
contradicts the connectivity of $\Gamma .$ $\Box $

\vskip 0.2truecm

The arguments contained in the proof of the next proposition will be used
many times in the rest of the paper.

\begin{proposition}
\label{prop1}: Let $\Gamma \subset D$ be a connected component and let $%
a\leq 0$ be such that $V_a$ intersects $\Gamma .$ Then, $\Gamma \cap V_a^{-}$
has at least one unlimited connected component, which intersects $V_a.$
\end{proposition}

{\it Proof:}

Let us consider the $L,R$-compactification of $\widetilde{A}={\rm I}%
\negthinspace {\rm R}\times [0,1],$ denoted $\widehat{A}.$ For every object
(point, set, etc) in $\widetilde{A},$ we denote the corresponding object in $
\widehat{A}$ by putting a $\widehat{}$ on it.

Let $z_n\in \Gamma \cap V_a^{-}$ be a sequence such that $p_1(z_n)\stackrel{%
n\rightarrow \infty }{\rightarrow }-\infty ,$ or equivalently, $\widehat{A}%
\ni \widehat{z}_n\stackrel{n\rightarrow \infty }{\rightarrow }L.$

Also, note that $\widehat{\Gamma }$ is connected, it intersects $\widehat{V}%
_a$ and contains $L.$ Each $\widehat{z}_n$ belongs to a connected component
of $\widehat{\Gamma }\cap \widehat{V}_a^{-},$ denoted $\widehat{\Gamma }_n.$
The connectivity of $\widehat{\Gamma }$ implies that each $\widehat{\Gamma }%
_n$ intersects $\widehat{V}_a.$ Let $\widehat{\Gamma }_{n_i}$ be a
convergent subsequence in the Hausdorff topology, $\widehat{\Gamma }_{n_i}%
\stackrel{n\rightarrow \infty }{\rightarrow }\widehat{\Gamma }^{*}.$ This
means that, given any open neighborhood $\widehat{N}$ of $\widehat{\Gamma }%
^{*},$ for all sufficiently large $i,$ $\widehat{\Gamma }_{n_i}$ is
contained in $\widehat{N}.$ So $\widehat{\Gamma }^{*}$ must contain $L$ and
must intersect $\widehat{V}_a.$ Suppose that $\widehat{\Gamma }^{*}$ is not
contained in $\widehat{\Gamma }.$ This means that there exists $\widehat{P}%
\in \widehat{\Gamma }^{*},$ with $\widehat{P}\notin \widehat{\Gamma }.$ As $
\widehat{\Gamma }$ is closed, for some $\epsilon _0>0,$ $B_{\epsilon _0}( 
\widehat{P})\cap \widehat{\Gamma }=\emptyset ,$ where $B_{\epsilon _0}( 
\widehat{P})=\{\widehat{z}\in \widehat{A}:d_{Euclidean}(\widehat{z},\widehat{%
P})<\epsilon _0\}$ and $d_{Euclidean}(\bullet ,\bullet )$ is the usual
Euclidean distance in $\widehat{A}.$ But as $\widehat{\Gamma }_{n_i}%
\stackrel{n\rightarrow \infty }{\rightarrow }\widehat{\Gamma }^{*}$ in the
Hausdorff topology, for all sufficiently large $i,$ $\widehat{\Gamma }%
^{*}\subset (\epsilon _0/2-neighborhood$ $of$ $\widehat{\Gamma }_{n_i}).$
Thus we get that $d_{Euclidean}(\widehat{P},\widehat{\Gamma })\leq
d_{Euclidean}(\widehat{P},\widehat{\Gamma }_{n_i})<\epsilon _0/2,$ something
that contradicts the choice of $\widehat{P}\in \widehat{\Gamma }^{*}.$ So $
\widehat{\Gamma }^{*}\subset \widehat{\Gamma }$ and the proposition is
proved because although $\Gamma ^{*}$ may not be connected, it must contain
an unlimited connected component which intersects $V_a$. $\Box $

\vskip 0.2truecm

Before going on, let us define the sets $\Gamma _{a,down}$ and $\Gamma
_{a,up}$ as follows:%
$$
\begin{array}{c}
\Gamma _{a,down}(\Gamma _{a,up})=\Gamma _{a,down}^{\prime }(\Gamma
_{a,up}^{\prime }) 
\text{ plus all points} \\ \text{ in the boundary of }\Gamma
_{a,down}^{\prime }(\Gamma _{a,up}^{\prime })\text{ of the form }(a, 
\widetilde{y}) 
\end{array}
$$
If $\Gamma $ is a connected component of $D$ which intersects some vertical $%
V_a,$ it is possible that $\Gamma \cap V_a^{-}$ has more then one unlimited
connected component. We denote by 
\begin{equation}
\label{defcolchete}\left[ \Gamma \cap V_a^{-}\right] =\text{union of all
unbounded connected components of }\Gamma \cap V_a^{-}. 
\end{equation}

\begin{proposition}
\label{prop2}:\ Let $a,b\in {\rm I}\negthinspace {\rm R}$ be such that $b<a$
and let $\Gamma $ be a connected component of $D,$ which intersects $V_a.$
Then $\Gamma _{b,down}\subset \Gamma _{a,down}$ and $\Gamma _{b,up}\subset
\Gamma _{a,up}.$
\end{proposition}

{\it Proof:}

Let $z\in \Gamma _{b,down}.$ This means that there exists a simple
continuous arc $\theta $ which connects $z$ to a point $z_0\in ]-\infty
,b[\times \{0\},${\bf \ $\theta \cap \Gamma =\emptyset $} and $\theta
\subset \Gamma _{b,down}\subset ]-\infty ,b[\times [0,1].$ As $a>b,$ $\theta
\cap \partial \Gamma _{a,down}\subset \theta \cap \Gamma =\emptyset .$ As $%
z_0\in \Gamma _{a,down},$ we get that $\theta \subset \Gamma _{a,down},$
which implies that $\Gamma _{b,down}\subset \Gamma _{a,down}.$ The other
inclusion is proved in a similar way. $\Box $

\vskip 0.2truecm

Let $\Gamma _1,\Gamma _2$ be two different connected components of $D$ and
let $V_a$ be a vertical which intersects $\Gamma _1.$

\begin{lemma}
\label{partecha2}: One and only one of the following possibilities must hold:

$\left[ \Gamma _2\cap V_a^{-}\right] \subset \Gamma _{1a,down}$ or $\left[
\Gamma _2\cap V_a^{-}\right] \subset \Gamma _{1a,up}.$
\end{lemma}

{\it Proof:}

First we prove that $\left[ \Gamma _2\cap V_a^{-}\right] \subset \Gamma
_{1a,down}\cup \Gamma _{1a,up}.$ Suppose this is not the case. Then, there
exists an unlimited connected component of $\Gamma _2\cap V_a^{-},$ denoted $%
\Gamma _2^{*},$ contained in a connected component of $\Gamma _1^{comp,a}$
(see (\ref{defgammacompa})) different from $\Gamma _{1a,down}$ and $\Gamma $$%
_{1a,up}.$ Denote this component by $\Gamma $$_{1a,mid}.$ Fix some $P\in
\Gamma _2^{*}.$ As $P\notin \Gamma _1,$ there exists $\epsilon >0$ such that 
$B_\epsilon (P)\cap \Gamma _1=\emptyset .$ Now, let $\alpha ^{\prime
}\subset {\rm I}\negthinspace {\rm R}\times ]0,1[$ be a simple continuous
arc which connects $P$ to $(1,0.5),$ totally contained in $\Gamma _1^c,$
which is an open connected set that contains $P$ and $(1,0.5).$ Moreover, as 
$]0,+\infty [\times [0,1]\subset \Gamma _1^c,$ we can take $\alpha ^{\prime
} $ so that it does not intersect $]1,+\infty [\times \{0.5\}.$ Finally, let 
$\alpha $ be a simple continuous arc given by $[1,+\infty [\times \{0.5\}$
plus a continuous part of $\alpha ^{\prime },$ whose endpoints are $(1,0.5)$
and some point in $\Gamma $$_2^{*},$ so that $\alpha \cap \Gamma $$_2^{*}$
consists of only its end point (clearly, this end point may not be $P$).

{\bf Properties of $\alpha \cup \Gamma_2^{*}:$}

\begin{itemize}
\item  {\bf $\alpha \cup \Gamma _2^{*}$} is a closed, connected set,
disjoint from ${\rm I}\negthinspace {\rm R\times \{0,1\};}$

\item  ${\rm I}\negthinspace {\rm R\times \{0\}}$ and ${\rm I}%
\negthinspace 
{\rm R\times \{1\}}$ are in different connected components of $\left( \alpha
\cup \Gamma _2^{*}\right) ^c;$

\item  $\alpha $ is limited to the left, that is, there exists a number $M>0$
such that, for all points $\widetilde{z}$ in $\alpha ,$ $p_1(\widetilde{z}%
)>-M;$

\item  $\left( \alpha \cup \Gamma _2^{*}\right) \cap ${\bf $\Gamma
_1=\emptyset ;$}
\end{itemize}

Let us choose $b<a$ such that $\alpha \subset V_{b+1/2}^{+}.$ By proposition 
\ref{prop2}, $\Gamma _{1b,down}\subset \Gamma _{1a,down}$ and $\Gamma
_{1b,up}\subset \Gamma _{1a,up},$ so we get that $\Gamma _2^{*}\cap \left(
\Gamma _{1b,down}\cup \Gamma _{1b,up}\right) =\emptyset .$ Now let $\beta
_0\subset \Gamma _{1b,down}$ and $\beta _1\subset \Gamma _{1b,up}$ be simple
continuous arcs which satisfy the following:

\begin{itemize}
\item  $\beta _0$ connects a point of $]-\infty ,b[\times \{0\}$ to a point
of $\Gamma _1;$

\item  $\beta _1$ connects a point of $]-\infty ,b[\times \{1\}$ to a point
of $\Gamma _1;$
\end{itemize}

So the following conditions hold%
$$
\left( \beta _0\cup \beta _1\right) \cap \Gamma _2^{*}=\emptyset \text{ and }%
\left( \beta _0\cup \beta _1\right) \subset V_b^{-}\text{ }\Rightarrow \text{
}\left( \beta _0\cup \beta _1\right) \cap \alpha =\emptyset 
$$

and thus 
$$
\left( \beta _0\cup \Gamma _1\cup \beta _1\right) \cap \left( \alpha \cup
\Gamma _2^{*}\right) =\emptyset , 
$$
something that contradicts the fact that $\left( \beta _0\cup \Gamma _1\cup
\beta _1\right) $ is a closed connected set and the\ ``{\bf Properties of $%
\alpha \cup \Gamma _2^{*}$}'' listed above, see figure 3. 
So, $\left[ \Gamma _2\cap
V_a^{-}\right] \subset (\Gamma _{1a,down}\cup \Gamma _{1a,up}).$

Suppose now that for some $\Gamma _2^{*},\Gamma _2^{**}\in \left[ \Gamma
_2\cap V_a^{-}\right] ,$ we have $\Gamma _2^{*}\subset \Gamma _{1a,down}$
and $\Gamma _2^{**}\subset \Gamma _{1a,up}.$ In the same way as above, there
exists a simple continuous arc $\alpha \subset {\rm I}\negthinspace {\rm R}%
\times ]0,1[$ which contains $[1,+\infty [\times \{0.5\}$ and connects some
point of $\Gamma _1$ to $(1,0.5),$ in a way that $\alpha \subset \Gamma _2^c$
and $\alpha $ intersects $\Gamma _1$ only at its end point$.$ Clearly, $%
\left( \alpha \cup \Gamma _1\right) $ is a closed connected set, which
satisfies:\ ${\rm I}\negthinspace {\rm R\times \{0\}}$ and ${\rm I}%
\negthinspace {\rm R\times \{1\}}$ are in different connected components of $%
\left( \alpha \cup \Gamma _1\right) ^c.$

Again, as above, let us choose $b<a$ such that $\alpha \subset V_{b+1}^{+}. $
Proposition \ref{prop1} implies that $\left[ \Gamma _2^{*}\cap
V_b^{-}\right] $ and $\left[ \Gamma _2^{**}\cap V_b^{-}\right] $ are
non-empty. From what we did above, we get that $\left[ \Gamma _2^{*}\cap
V_b^{-}\right] \cup \left[ \Gamma _2^{**}\cap V_b^{-}\right] \subset \Gamma
_{1b,down}\cup \Gamma _{1b,up}.$

If, $\left[ \Gamma _2^{*}\cap V_b^{-}\right] \cap \Gamma _{1b,up}\neq
\emptyset $ $\Rightarrow $ $\Gamma _2^{*}\cap \Gamma _{1b,up}\neq \emptyset
, $ which implies, by proposition \ref{prop2}, that $\Gamma _2^{*}\cap
\Gamma _{1a,up}\neq \emptyset ,$ a contradiction. So, $\left[ \Gamma
_2^{*}\cap V_b^{-}\right] \cap \Gamma _{1b,up}=\emptyset $ and a similar
argument gives $\left[ \Gamma _2^{**}\cap V_b^{-}\right] \cap \Gamma
_{1b,down}=\emptyset .$ So, $\left[ \Gamma _2^{*}\cap V_b^{-}\right] \subset
\Gamma _{1b,down}$ and $\left[ \Gamma _2^{**}\cap V_b^{-}\right] \subset
\Gamma _{1b,up}.$ Thus, there exists a simple continuous arc $\beta _0$
contained in $\Gamma _{1b,down}$ which connects a point of $\Gamma _2^{*}$
to some point in $]-\infty ,b[\times \{0\}.$ Similarly, there exists a
simple continuous arc $\beta_1$ contained in $\Gamma _{1b,up}$ which
connects a point of $\Gamma _2^{**}$ to some point in $]-\infty ,b[\times
\{1\},$ see figure 4. But 
$\left( \beta _0\cup \Gamma _2\cup \beta _1\right) $ is a closed
connected set and by construction of $\beta _0$ and $\beta _1,$%
$$
\left( \beta _0\cup \Gamma _2\cup \beta _1\right) \cap \left( \alpha \cup
\Gamma _1\right) =\emptyset , 
$$
which is a contradiction, completing the proof of the lemma. $\Box $

\vskip 0.2truecm

The previous results will be used in what follows in order to define a
complete ordering among the connected components of $D.$

Let $\Gamma _1$ and $\Gamma _2$ be different connected components of $D$ and
let $a\in {\rm I}\negthinspace {\rm R}$ be such that $\Gamma _1$ and $\Gamma
_2$ intersect $V_a.$ We say that $\Gamma _2\prec _a\Gamma _1,$ if $\left[
\Gamma _2\cap V_a^{-}\right] \subset \Gamma _{1a,down}.$

\begin{lemma}
\label{partecha3}: Given $\Gamma _1,\Gamma _2$ and $a\in {\rm I}%
\negthinspace {\rm R}$ as above, either $\Gamma _2\prec _a\Gamma _1$ or $%
\Gamma _1\prec _a\Gamma _2.$
\end{lemma}

{\it Proof: }

From lemma \ref{partecha2}, either $\left[ \Gamma _2\cap V_a^{-}\right]
\subset \Gamma _{1a,down}$ or $\left[ \Gamma _2\cap V_a^{-}\right] \subset
\Gamma _{1a,up.}$ In the first possibility, $\Gamma _2\prec _a\Gamma _1.$ So
we are left to show that, if $\left[ \Gamma _2\cap V_a^{-}\right] \subset
\Gamma _{1a,up},$ then $\left[ \Gamma _1\cap V_a^{-}\right] \subset \Gamma
_{2a,down},$ which means that $\Gamma _1\prec _a\Gamma _2.$

Thus, let us suppose that $\left[ \Gamma _2\cap V_a^{-}\right] \subset
\Gamma _{1a,up}$ and $\left[ \Gamma _1\cap V_a^{-}\right] \subset \Gamma
_{2a,up}.$ If we arrive at a contradiction, the lemma will be proved.

The argument here is very similar to the one used in the proof of lemma \ref
{partecha2}. First, choose a simple continuous arc $\alpha \subset {\rm I}%
\negthinspace {\rm R}\times ]0,1[$ which contains $[1,+\infty [\times
\{0.5\} $ and connects some point of $\Gamma _1$ to $(1,0.5),$ in a way that 
$\alpha \subset \Gamma _2^c$ and $\alpha $ intersects $\Gamma _1$ only at
its end point$.$ Clearly, $\left( \alpha \cup \Gamma _1\right) $ is a closed
connected set, which satisfies:\ ${\rm I}\negthinspace {\rm R\times \{0\}}$
and ${\rm I}\negthinspace {\rm R\times \{1\}}$ are in different connected
components of $\left( \alpha \cup \Gamma _1\right) ^c.$

As $\left[ \Gamma _2\cap V_a^{-}\right] \subset \Gamma _{1a,up},$ there
exists an element of $\left[ \Gamma _2\cap V_a^{-}\right] ,$ denoted $\Gamma
_2^{*},$ which by definition is closed, connected, unlimited to the left and
is contained in $\Gamma _{1a,up}.$ Again, let us choose $b<a$ such that $%
\alpha \subset V_{b+1}^{+}.$

As in the end of the proof of lemma \ref{partecha2}, we get that $\left[
\Gamma _2^{*}\cap V_b^{-}\right] \subset \Gamma _{1b,up}.$ So, there exists
a simple continuous arc $\beta_1$ contained in $\Gamma _{1b,up}$ which
connects a point of $\Gamma _2^{*}$ to some point in $]-\infty ,b[\times
\{1\}.$ Clearly, $\beta _1\cap \Gamma _1=\emptyset .$

As $\left[ \Gamma _1\cap V_a^{-}\right] \subset \Gamma _{2a,up},$ an
argument similar to the one used to prove proposition \ref{prop1} implies:

\begin{proposition}
\label{limitado}: There exists a real number $c\leq b,$ such that $\left(
\Gamma _1\cap V_c^{-}\right) \cap \Gamma _{2a,down}=\emptyset .$
\end{proposition}

{\it Proof:}

Suppose by contradiction, that the fact is not true. Then, there is a
sequence of points $z_n\in \Gamma _1\cap \Gamma _{2a,down},$ such that $%
p_1(z_n)\stackrel{n\rightarrow \infty }{\rightarrow }-\infty ,$ or
equivalently, $\widehat{A}\ni \widehat{z}_n\stackrel{n\rightarrow \infty }{%
\rightarrow }L.$ Each $\widehat{z}_n$ belongs to a connected component of $
\widehat{\Gamma }_1\cap \widehat{V}_a^{-},$ denoted $\widehat{\Gamma }%
_{1,n}\subset closure\left( \widehat{\Gamma }_{2a,down}\right) .$ The
connectivity of $\Gamma _1$ and the fact that it intersects $V_a$ implies
that each $\widehat{\Gamma }_{1,n}$ intersects $\widehat{V}_a.$ Let $
\widehat{\Gamma }_{1,n_i}$ be a convergent subsequence in the Hausdorff
topology, $\widehat{\Gamma }_{1,n_i}\stackrel{n\rightarrow \infty }{%
\rightarrow }\widehat{\Gamma }^{*}.$ This means that, given any open
neighborhood $\widehat{N}$ of $\widehat{\Gamma }^{*},$ for all sufficiently
large $i,$ $\widehat{\Gamma }_{1,n_i}$ is contained in $\widehat{N}.$ So $
\widehat{\Gamma }^{*}$ must contain $L$ and must intersect $\widehat{V}_a.$
Suppose that $\widehat{\Gamma }^{*}$ is not contained in $\widehat{\Gamma }%
_1.$ This means that there exists $\widehat{P}\in \widehat{\Gamma }^{*},$
with $\widehat{P}\notin \widehat{\Gamma }_1.$ As $\widehat{\Gamma }_1$ is
closed, for some $\epsilon _0>0,$ $B_{\epsilon _0}(\widehat{P})\cap \widehat{%
\Gamma }_1=\emptyset ,$ where $B_{\epsilon _0}(\widehat{P})=\{\widehat{z}\in 
\widehat{A}:d_{Euclidean}(\widehat{z},\widehat{P})<\epsilon _0\}$ and $%
d_{Euclidean}(\bullet ,\bullet )$ is the usual Euclidean distance in $
\widehat{A}.$ But as $\widehat{\Gamma }_{n_i}\stackrel{n\rightarrow \infty }{%
\rightarrow }\widehat{\Gamma }^{*}$ in the Hausdorff topology, for all
sufficiently large $i,$ $\widehat{\Gamma }^{*}\subset (\epsilon
_0/2-neighborhood$ $of$ $\widehat{\Gamma }_{n_i}).$ Thus we get that $%
d_{Euclidean}(\widehat{P},\widehat{\Gamma }_1)\leq d_{Euclidean}(\widehat{P}%
, \widehat{\Gamma }_{1,n_i})<\epsilon _0/2,$ something that contradicts the
choice of $\widehat{P}\in \widehat{\Gamma }^{*}.$ So $\widehat{\Gamma }%
^{*}\subset \widehat{\Gamma }_1.$ Clearly $\widehat{\Gamma }^{*}\subset
closure\left( \widehat{\Gamma }_{2a,down}\right) $ and, as $\Gamma _1\cap
\Gamma _2=\emptyset ,$ we get, by lemma \ref{partecha2} that $\left[ \Gamma
_1\cap V_a^{-}\right] \subset \Gamma _{2a,down},$ a contradiction with our
hypothesis. $\Box $

\vskip 0.2truecm

Now let us look at $\Gamma _{2c,down}\subset \Gamma _{2b,down},$ where $c$
comes from proposition \ref{limitado}. Thus, $\Gamma _1\cap \Gamma
_{2c,down}=\emptyset .$ So, there exists a simple continuous arc $\beta _0$
which connects a point of $\Gamma _2$ to some point in $]-\infty ,c[\times
\{0\},$ in a way that $\beta _0\cap \Gamma _2$ is one extreme of $\beta _0,$
denoted $m_0,$ and $\beta _0\backslash \{m_0\}\subset \Gamma _{2c,down},$
which implies that $\beta _0\cap \Gamma _1=\emptyset $ and $\beta _0\cap
\alpha =\emptyset .$ So, $\left( \beta _0\cup \Gamma _2\cup \beta _1\right) $
is a closed connected set, which intersects ${\rm I}\negthinspace {\rm %
R\times \{0\}}$ and ${\rm I}\negthinspace {\rm R\times \{1\}.}$ And, by
construction

$$
\left( \beta _0\cup \Gamma _2\cup \beta _1\right) \cap \left( \alpha \cup
\Gamma _1\right) =\emptyset , 
$$
a contradiction. So if $\left[ \Gamma _2\cap V_a^{-}\right] \subset \Gamma
_{1a,up},$ then $\left[ \Gamma _1\cap V_a^{-}\right] \subset \Gamma
_{2a,down},$ which implies that $\Gamma _1\prec _a\Gamma _2$ and the lemma
is proved. $\Box $

\vskip 0.2truecm

Finally, in order to present a good definition of order, we need the
following lemma:

\begin{lemma}
\label{ordem}: Let $\Gamma _1$ and $\Gamma _2$ be different connected
components of $D$ and let $a,b\in {\rm I}\negthinspace {\rm R}$ be such that 
$\Gamma _1$ and $\Gamma _2$ intersect $V_a$ and $V_b.$ Then we have the
following:%
$$
\begin{array}{c}
\Gamma _1\prec _a\Gamma _2
\text{ }\Leftrightarrow \text{ }\Gamma _1\prec _b\Gamma _2 \\ \Gamma _2\prec
_a\Gamma _1\text{ }\Leftrightarrow \text{ }\Gamma _2\prec _b\Gamma _1
\end{array}
$$
\end{lemma}

{\it Proof:}

Suppose that $b<a$ and $\Gamma _2\prec _a\Gamma _1$ $\Leftrightarrow $ $%
\left[ \Gamma _2\cap V_a^{-}\right] \subset \Gamma _{1a,down}.$ Proposition 
\ref{limitado} tells us that $\left( \Gamma _2\cap V_a^{-}\right) \cap
\Gamma _{1a,up}$ is a limited set. So, as $\Gamma _{1b,up}\subset \Gamma
_{1a,up},$ $\left[ \Gamma _2\cap V_b^{-}\right] $ must be contained in $%
\Gamma _{1b,down},$ which means that $\Gamma _2\prec _b\Gamma _1.$ The other
implications are proved in a similar way. $\Box $

\vskip 0.2truecm

So given $\Gamma _1$ and $\Gamma _2,$ two different connected components of $%
D,$ if $a\in {\rm I}\negthinspace {\rm R}$ is such that $\Gamma _1$ and $%
\Gamma _2$ intersect $V_a,$ we can define an order $\prec $ between them as
explained above and this order is independent of the choice of $a.$ We only
need the following condition: $V_a$ must intersect $\Gamma _1$ and $\Gamma
_2.$

Also, let us prove the following associativity lemma:

\begin{lemma}
\label{associati}:\ If $\Gamma _1,\Gamma _2,\Gamma _3$ are connected
components of $D,$ such that $\Gamma _1\prec \Gamma _2$ and $\Gamma _2\prec
\Gamma _3,$ then $\Gamma _1\prec \Gamma _3.$
\end{lemma}

{\it Proof:}

Let $a\in {\rm I}\negthinspace {\rm R}$ be such that $\Gamma _1,\Gamma _2$
and $\Gamma _3$ intersect $V_a.$ Then, $\left[ \Gamma _1\cap V_a^{-}\right]
\subset \Gamma _{2a,down}$ and $\left[ \Gamma _2\cap V_a^{-}\right] \subset
\Gamma _{3a,down}.$ In the proof of lemma \ref{partecha3}, we proved that if 
$\Theta $ and $\Lambda $ are different connected components of $D$ and $a\in 
{\rm I}\negthinspace {\rm R}$ is such that $\Theta $ and $\Lambda $
intersect $V_a$ then, $\left[ \Theta \cap V_a^{-}\right] \subset \Lambda
_{a,down}$ implies $\left[ \Lambda \cap V_a^{-}\right] \subset \Theta
_{a,up}.$ So, $\left[ \Gamma _2\cap V_a^{-}\right] \subset \Gamma _{1a,up}$
and $\left[ \Gamma _3\cap V_a^{-}\right] \subset \Gamma _{2a,up}.$ Now,
using proposition \ref{limitado}, let us choose $b\leq a$ such that the
following inclusions hold: 
\begin{equation}
\label{cacaca} 
\begin{array}{l}
\Gamma _3\cap V_b^{-}\subset \Gamma _{2b,up} \\ 
\Gamma _1\cap V_b^{-}\subset \Gamma _{2b,down} \\ 
\Gamma _2\cap V_b^{-}\subset \Gamma _{3b,down} 
\end{array}
\end{equation}

Finally, let us prove that $\Gamma _{2b,down}\subset \Gamma _{3b,down}.$

If this is not the case, then there exists a simple continuous arc $\alpha
\subset \Gamma _{2b,down}$ that connects a point from $]-\infty ,b[\times
\{0\}$ to a point $P\notin \Gamma _{3b,down}.$ Thus $\alpha $ intersects $%
\Gamma _3,$ a contradiction with expression (\ref{cacaca}). So, $\Gamma
_1\cap V_b^{-}\subset \Gamma _{2b,down}\subset \Gamma _{3b,down},$ which
implies that $\left[ \Gamma _1\cap V_b^{-}\right] \subset \Gamma _{3b,down}$ 
$\Leftrightarrow $ $\Gamma _1\prec \Gamma _3$ and the lemma is proved. $\Box 
$

\vskip 0.2truecm

Our next objective is to show that $\widetilde{f}$ preserves the order just
defined. First, note that if $\Gamma $ is a connected component of $D,$ as $
\widetilde{f}(D)\subset D$ and $\widetilde{f}(\Gamma )$ is connected,
unlimited to the left, there exists a connected component of $D,$ denoted $%
\Gamma ^{+}$, that contains $\widetilde{f}(\Gamma ).$

We have

\begin{lemma}
\label{presordem}: Let $\Gamma _1,\Gamma _2$ be connected components of $%
B^{-}$ and suppose $\Gamma _1\prec \Gamma _2.$ Then, $\widetilde{f}(\Gamma
_1)\prec \widetilde{f}(\Gamma _2)$ and so, if $\Gamma _1^{+}\neq \Gamma
_2^{+},$ then $\Gamma _1^{+}\prec \Gamma _2^{+}.$
\end{lemma}

{\it Proof:}

Suppose that $\Gamma _2^{+}\prec \Gamma _1^{+}.$ As $\Gamma _1\prec \Gamma
_2,$ for any $a\in {\rm I}\negthinspace {\rm R}$ such that $\Gamma _1$ and $%
\Gamma _2$ intersect $V_a,$ the proof of lemma \ref{partecha3} implies that $%
\left[ \Gamma _2\cap V_a^{-}\right] \subset \Gamma _{1a,up}.$ From
proposition \ref{limitado}, there exists a sufficiently small $b<0$ such
that: 
\begin{equation}
\label{estcont} 
\begin{array}{l}
\Gamma _1\cap V_b^{-}\subset \Gamma _{2b,down} \\ 
\Gamma _2\cap V_b^{-}\subset \Gamma _{1b,up} \\ 
\Gamma _2^{+}\cap V_b^{-}\subset \Gamma _{1b,down}^{+} 
\end{array}
\end{equation}

Let $c<b$ be such that $\widetilde{f}^{\pm 1}(V_c)\cap V_b=\emptyset .$ From
our previous results, we get that%
$$
\begin{array}{l}
\Gamma _1\cap V_c^{-}\subset \Gamma _{2c,down} \\ 
\Gamma _2\cap V_c^{-}\subset \Gamma _{1c,up} \\ 
\Gamma _2^{+}\cap V_c^{-}\subset \Gamma _{1c,down}^{+} 
\end{array}
. 
$$
So, there exists a simple continuous arc $\alpha \subset \Gamma
_{1c,down}^{+}\subset V_c^{-}$ that connects a point from $]-\infty
,c[\times \{0\}$ to a point $P\in \widetilde{f}(\Gamma _2)\subset \Gamma
_2^{+}.$ From the choice of $c,$ $\widetilde{f}^{-1}(\alpha )\subset V_b^{-}$
and it connects a point from $]-\infty ,b[\times \{0\}$ to $\widetilde{f}%
^{-1}(P)\in \Gamma _2.$ As $\alpha \subset \Gamma _{1c,down}^{+},$ $\alpha
\cap \Gamma _1^{+}=\emptyset ,$ so $\widetilde{f}^{-1}(\alpha )\cap \Gamma
_1=\emptyset .$ Thus, $\widetilde{f}^{-1}(\alpha )\subset \Gamma _{1b,down},$
which implies that $\Gamma _2\cap \Gamma _{1b,down}\neq \emptyset $ and this
contradicts (\ref{estcont}). So, either $\Gamma _1^{+}=\Gamma _2^{+},$ or $%
\Gamma _1^{+}\prec \Gamma _2^{+},$ because what the proof presented above
really shows is that

$$
\Gamma _1\prec \Gamma _2\Rightarrow \widetilde{f}(\Gamma _1)\prec \widetilde{%
f}(\Gamma _2).\text{ }\Box 
$$

\vskip 0.2truecm

In the particular cases where $D=B^{-}$ or when $D=\omega (B^{-}),$ we can
show that $\widetilde{f}^{-1}$ is also an order-preserving transformation.
This is more clearly seen when $D=\omega (B^{-})$ since $\omega (B^{-})$ is $
\widetilde{f}$-invariant. If $\Gamma $ is a connected component of $\omega
(B^{-})$, so is $\widetilde{f}^{-1}(\Gamma ),$ which we call $\Gamma ^{-}.$
It is then a simple consequence of the previous lemma that, if $\Gamma
_1,\Gamma _2\subset \omega (B^{-}),$ with $\Gamma _1\prec \Gamma _2,$ then $%
\Gamma _1^{-}\prec \Gamma _2^{-}.$

When $D=B^{-},$ if $\Gamma $ is a connected component of $B^{-},$ then $
\widetilde{f}^{-1}(\Gamma )$ is also closed connected, unlimited to the
left, and limited to the right. The only possibility that may prevent $
\widetilde{f}^{-1}(\Gamma )$ from being contained in a connected component
of $B^{-}$ is the following: $\widetilde{f}^{-1}(\Gamma )\cap ]0,+\infty
[\times [0,1]\neq \emptyset $. There are 2 possibilities:

\begin{enumerate}
\item  $\widetilde{f}^{-1}(\Gamma )\cap ]0,+\infty [\times [0,1]=\emptyset .$
As $\widetilde{f}\left( \widetilde{f}^{-1}(\Gamma )\right) =\Gamma \subset
B^{-},$ there is a connected component of $B^{-},$ denoted $\Gamma ^{-},$
which satisfies $\Gamma ^{-}\supset \widetilde{f}^{-1}(\Gamma )$ and $
\widetilde{f}\left( \Gamma ^{-}\right) \cap \Gamma \neq \emptyset .$ This
implies that$\Rightarrow $ $\widetilde{f}\left( \Gamma ^{-}\right) \subset
\Gamma ,$ and so $\Gamma ^{-}=\widetilde{f}^{-1}(\Gamma ).$

\item  $\widetilde{f}^{-1}(\Gamma )\cap ]0,+\infty [\times [0,1]\neq
\emptyset .$ As $\left[ \widetilde{f}^{-1}(\Gamma )\cap V^{-}\right] $ has
at least one connected component, denoted $\Gamma ^{*},$ we get that $
\widetilde{f}\left( \Gamma ^{*}\right) \subset \Gamma \subset B^{-},$ so
there is a connected component of $B^{-},$ denoted $\Gamma ^{-},$ which
contains $\Gamma ^{*}$ and $\widetilde{f}\left( \Gamma ^{-}\right) \subset
\Gamma $ because $\widetilde{f}\left( \Gamma ^{-}\right) \cap \Gamma \neq
\emptyset .$ Note that in this case, $\Gamma ^{-}$ may not be unique.
\end{enumerate}

We can still formulate the following result:

\begin{lemma}
\label{presordem2}: Let $\Gamma _1,\Gamma _2$ be connected components of $D$
and suppose $\Gamma _1\prec \Gamma _2.$ Then, for any choice of $\Gamma
_1^{-}$ and $\Gamma _2^{-}$, we have $\Gamma _1^{-}\prec \Gamma _2^{-}.$
\end{lemma}

{\it Proof:} Using lemma \ref{presordem}, as $\Gamma _1\neq \Gamma _2,$ if $%
\Gamma _2^{-}\prec \Gamma _1^{-},$ then $\Gamma _2\prec \Gamma _1,$
therefore we must have $\Gamma _1^{-}\prec \Gamma _2^{-}.$ $\Box $

\vskip 0.2truecm

We now, for a fixed connected component $\Gamma $ of $D,$ consider the
covering mapping $p\mid _\Gamma $. It may or may not be injective. We
examine the consequences in each case:

\subsubsection{The covering mapping $p\mid _\Gamma $ is not injective}

This means that there exists $\widetilde{z}\in \widetilde{A}$ and an integer 
$s>0$ such that $\widetilde{z},\widetilde{z}+(s,0)\in \Gamma .$ So, $\Gamma
\cap \Gamma -(s,0)\neq \emptyset .$ The last property of $D$ tell us that $%
\Gamma -(s,0)\subset D.$ But this implies that 
\begin{equation}
\label{gamanonin}\Gamma -(s,0)\subset \Gamma , 
\end{equation}
because $\Gamma $ is a connected component of $D.$

Suppose that $\Gamma -(1,0)$ is not contained in $\Gamma .$ As $\Gamma
-(1,0)\subset D,$ we get that $\Gamma -(1,0)\cap \Gamma =\emptyset .$ As $%
\Gamma -(1,0)$ does not intersect $V_{m_\Gamma }=\{m_\Gamma \}\times [0,1],$
lemma \ref{partecha2} implies that either $\Gamma -(1,0)\subset \Gamma
_{down}$ or $\Gamma -(1,0)\subset \Gamma _{up}.$ Suppose it is contained in $%
\Gamma _{up}.$

\begin{proposition}
\label{gamaup}:\ If $\Gamma -(1,0)\subset \Gamma _{up},$ then $\Gamma
-(i,0)\subset \Gamma _{up}$ for all integers $i>1.$
\end{proposition}

{\it Proof:}

Suppose there exists $s_0>1$ (the smallest one) such that $\Gamma
-(s_0,0)\subset \Gamma .$ This means that $\Gamma ,\Gamma -(1,0),...,\Gamma
-(s_0-1,0)$ are all disjoint.

As $\Gamma -(1,0)\subset \Gamma _{up}$ (which implies that $\Gamma \prec $ $%
\Gamma -(1,0)),$ we get that $\Gamma -(s,0)\cap $ $\Gamma -(s+1,0)=\emptyset 
$ and $\Gamma -(s,0)\prec $ $\Gamma -(s+1,0),$ for all integers $s>0.$ So,
in particular, using lemma \ref{associati}, we obtain the following
implications: 
\begin{equation}
\label{expnova} 
\begin{array}{l}
1) 
\text{ }\Gamma \prec \Gamma -(1,0)\prec \Gamma -(2,0)\prec \Gamma
-(3,0)\prec ...\prec \Gamma -(s_0-1,0) \\ 2)\text{ }\Gamma -(s_0-1,0)\prec
\Gamma -(s_0,0). 
\end{array}
\end{equation}

So, as $\Gamma -(s_0,0)\subset \Gamma $ and $\Gamma \cap $ $\Gamma
-(s_0-1,0)=\emptyset ,$ we get from 2) of (\ref{expnova}) that $\Gamma
-(s_0-1,0)\prec \Gamma ,$ a contradiction with 1) of (\ref{expnova}). Thus
for all integers $\dot i>0,$ $\Gamma \cap $ $\Gamma -(i,0)=\emptyset $ and
so 
$$
\Gamma \prec \Gamma -(1,0)\prec \Gamma -(2,0)\prec \Gamma -(3,0)\prec
...\prec \Gamma -(i,0). 
$$
If $a\in {\rm I}\negthinspace {\rm R}$ is such that $\Gamma $ and $\Gamma
-(i,0)$ intersect $V_a,$ then, $\left[ \Gamma \cap V_a^{-}\right] \subset
\Gamma -(i,0)_{a,down}$ $\Leftrightarrow $ $\left[ \Gamma -(i,0)\cap
V_a^{-}\right] \subset \Gamma _{a,up}.$ As $\Gamma _{a,up}\subset \Gamma
_{up},$ we get that $\Gamma -(i,0)\cap \Gamma _{up}\neq \emptyset ,$ and so $%
\Gamma -(i,0)\subset \Gamma _{up},$ because $\Gamma -(i,0)\cap V_{m_\Gamma
}=\emptyset .$ $\Box $

\vskip 0.2truecm

Therefore, if the map $p\mid_{\Gamma}$ is not injective, then 
\begin{equation}
\label{gamma-1_em_gamma} \Gamma - (1,0) \subset \Gamma. 
\end{equation}

We will call $\Gamma $ a non-injective component.

\subsubsection{The covering mapping $p\mid _\Gamma $ is injective}

This implies that $\Gamma \cap \Gamma +(s,0)=\emptyset ,$ for all integers $%
s\neq 0.$ In particular, $\Gamma \cap \Gamma -(1,0)=\emptyset $ and we use
this relation to describe the asymptotic behavior of $p(\Gamma )$ around the
annulus. As we explained just before defining the order $\prec $, any two
unlimited closed connected disjoint subsets of $V^{-}\subset \widetilde{A}$
which have connected complements, denoted $\Theta _1$ and $\Theta _2,$ are
related by $\prec ,$ that is, either $\Theta _1\prec \Theta _2$ or $\Theta
_2\prec \Theta _1.$ So, we say that $\Gamma $ is a down component of $D$ if $%
\Gamma \prec \Gamma -(1,0)$ and, analogously, $\Gamma $ is an up component
if $\Gamma -(1,0)$$\prec \Gamma .$

\begin{lemma}
\label{acumula}:\ If $\Gamma \subset D$ is a down component, then $%
dist(\Gamma ,$ ${\rm I}\negthinspace {\rm R}\times \{1\})>0$ and
analogously, \ if $\Gamma \subset D$ is an up component, then $dist(\Gamma ,$
${\rm I}\negthinspace {\rm R}\times \{0\})>0.$
\end{lemma}

{\it Proof:}

In both cases, the proof is analogous, so suppose $\Gamma $ is a down
component. This means that $\Gamma -(1,0)$ is contained in $\Gamma _{up}.$

Thus, for any $\widetilde{x}<m_\Gamma $ (see(\ref{mgama})), if we consider
the segment $\{\widetilde{x}\}\times [0,\widetilde{y}^{*}],$ where 
\begin{equation}
\label{ytilast}\widetilde{y}^{*}=\widetilde{y}^{*}(\widetilde{x})=\sup \{ 
\widetilde{y}\in ]0,1[:\Gamma \cap \{\widetilde{x}\}\times [0,\widetilde{y}%
]=\emptyset \}, 
\end{equation}
we get that $\Gamma -(1,0)\cap \{\widetilde{x}\}\times [0,\widetilde{y}%
^{*}]=\emptyset .$

Now, consider a point $(m_\Gamma -1,\widetilde{y}_\Gamma )\in \Gamma -(1,0)$
and a simple continuous arc $\gamma \subset int(\widetilde{A}),$ such that:

i) $\gamma \cap \left(\Gamma -(1,0)\right)=(m_\Gamma -1,\widetilde{y}_\Gamma
)$

ii) $\gamma \cap \Gamma =\emptyset $

iii) the endpoints of $\gamma $ are $(m_\Gamma -1,\widetilde{y}_\Gamma )$
and $(m_\Gamma +1,0.5)$

iv) $\gamma \cap \{m_\Gamma +1\}\times [0,1]=(m_\Gamma +1,0.5)$

As $\Gamma -(1,0)\subset \Gamma _{up}$ and $\left( \Gamma \cup \Gamma
-(1,0)\right) ^c$ is connected, it is possible to choose $\gamma $ as above,
see figure 5.

The complement of the closed connected set $\Gamma -(1,0)\cup \gamma \cup
\{m_\Gamma +1\}\times [0,1]$ has exactly two connected components in $%
]-\infty ,m_\Gamma +1[\times [0,1],$ one containing $]-\infty ,m_\Gamma
+1[\times \{0\},$ denoted $\Gamma -(1,0)_{down}$ and the other containing $%
]-\infty ,m_\Gamma +1[\times \{1\},$ denoted $\Gamma -(1,0)_{up}.$ Note that
this construction is not unique, because we may have more then one point in $%
\Gamma -(1,0)\cap \{m_\Gamma -1\}\times [0,1].$ Nevertheless, for any such
choice, $\Gamma \subset \Gamma -(1,0)_{down}.$ This follows from the fact
that, for any 
$$
\widetilde{x}<\min \left[ m_\Gamma ,\min \{\widetilde{x}\in {\rm I}%
\negthinspace {\rm R:}\ (\widetilde{x},\widetilde{y})\in \gamma ,\ \text{for
some}\ \widetilde{y}\in ]0,1[\}\right] -10, 
$$
the segment $\{\widetilde{x}\}\times [0,\widetilde{y}^{*}]$ (see (\ref
{ytilast})) does not intersect $\Gamma -(1,0)\cup \gamma \cup \{m_\Gamma
+1\}\times [0,1]$ and $(\widetilde{x},\widetilde{y}^{*})\in \Gamma .$


Now suppose, by contradiction, that $dist(\Gamma ,$ ${\rm I}\negthinspace 
{\rm R}\times \{1\})=0.$ As $\Gamma $ is closed and $\Gamma \cap {\rm I}%
\negthinspace {\rm R}\times \{1\}=\emptyset ,$ we get that for every 
$$
M\leq M_0=\min \left[ m_\Gamma -10,\min \{\widetilde{x}\in {\rm I}%
\negthinspace {\rm R:}\ (\widetilde{x},\widetilde{y})\in \gamma ,\ \text{for
some}\ \widetilde{y}\in ]0,1[\}\right] -10 
$$
there exists $\epsilon >0$ such that if $\widetilde{z}\in \Gamma \cap {\rm I}%
\negthinspace 
{\rm R\times }[1-\epsilon ,1],$ then $p_1(\widetilde{z})<M.$ So, for $M_0$
and $\epsilon >0$ as above, let us choose a point $\widetilde{z}_0\in \Gamma
\cap {\rm I}\negthinspace {\rm R\times }[1-\epsilon ,1]$ such that $p_1( 
\widetilde{z}_0)\geq p_1(\widetilde{z})$ for all $\widetilde{z}\in \Gamma
\cap {\rm I}\negthinspace {\rm R\times }[1-\epsilon ,1]$ and 
$$
\begin{array}{c}
dist( 
\widetilde{z}_0,{\rm I}\negthinspace {\rm R}\times \{1\})<dist(\widetilde{z},%
{\rm I}\negthinspace {\rm R}\times \{1\})\text{ for all } \\ \widetilde{z}%
\in \Gamma \cap \{p_1(\widetilde{z}_0)\}\times [1-\epsilon ,1]\text{ with } 
\widetilde{z}\neq \widetilde{z}_0. 
\end{array}
$$

Intuitively, if we start going left from $\{m_\Gamma \}\times [0,1], 
\widetilde{z}_0$ is the point of $\Gamma $ with largest possible $\widetilde{%
x}$ and $\widetilde{y}$ coordinates, that belongs to ${\rm I}\negthinspace 
{\rm R\times }[1-\epsilon ,1]$.

Now consider a closed vertical segment $l$ contained in ${\rm I}%
\negthinspace {\rm R\times }[1-\epsilon ,1],$ starting at $\widetilde{z}_0$
and ending at ${\rm I}\negthinspace {\rm R\times }\{1\}.$ By construction of 
$l,$ $l\cap \Gamma =\widetilde{z}_0.$ As $\Gamma \subset \Gamma
-(1,0)_{down} $ and $l\cap \left( \gamma \cup \{m_\Gamma +1\}\times
[0,1]\right) =\emptyset ,$ we get that $l\cap \Gamma -(1,0)\neq \emptyset .$
So, there exists $\widetilde{z}_1\in l\cap \Gamma -(1,0)$ which implies that 
$\widetilde{z}_1+(1,0)\in l+(1,0)\cap \Gamma .$ And this contradicts the
choice of $\widetilde{z}_0.$ $\Box $

\vskip 0.2truecm

In this case we will say that $\Gamma $ is an injective component.

\section{Proof of theorem 2}

Suppose, by contradiction, that $\omega (B^{-})\neq \emptyset .$ This
implies, by lemma \ref{conts1}, that either $S^1\times \{0\}\subset 
\overline{p(\omega (B^{-}))}$ or $S^1\times \{1\}\subset \overline{p(\omega
(B^{-}))}.$ Let us assume, without loss of generality, that $S^1\times
\{0\}\subset \overline{p(\omega (B^{-}))}.$

Since the rotation number of $\widetilde{f}$ restricted to $S^1\times \{0\}$
is strictly positive, there exists $\sigma >0$ such that $p_1(\widetilde{f}( 
\widetilde{x},0))>\widetilde{x}+2\sigma $ for all $\widetilde{x}\in {\rm I}%
\negthinspace 
{\rm R}.$ Let $\epsilon >0$ be sufficiently small such that for all $( 
\widetilde{x},\widetilde{y})\in {\rm I}\negthinspace {\rm R}\times
[0,\epsilon ],$ $p_1\circ \tilde f(\tilde x,\widetilde{y})>\tilde x+\sigma .$

As $S^1\times \{0\}\subset \overline{p(\omega(B^{-}))}$, there is a real $a$
such that 
\begin{equation}
\label{definia}\omega(B^{-}) \cap \{a\}\times [0,\epsilon ]\neq \emptyset . 
\end{equation}

The fact that $\omega (B^{-})$ is closed implies that there must be a $%
\delta \le \epsilon $ such that $(a,\delta )\in \omega (B^{-})$, and such
that for all $0\le \tilde y<\delta ,(a,\tilde y)\notin \omega (B^{-}).$ In
other words, $(a,\delta )$ is the ``lowest'' point of $\omega (B^{-})$ in $%
\{a\}\times [0,\epsilon ].$ We denote by $v$ the segment $\{a\}\times
[0,\delta [.$

Let $\Theta _1$ be the (unbounded) connected component of $\omega (B^{-})$
that contains $(a,\delta ).$ Let $\Omega $ be the connected component of $%
\left( \Theta _1\cup v\right) ^c$ that contains $]-\infty ,a[\times \{0\}.$
Of course, $\partial \Omega \subset \Theta _1\cup v,$ and $\partial 
\widetilde{f}(\Omega )\subset \widetilde{f}(\Theta _1)\cup \widetilde{f}(v).$

Note that, since $\omega (B^{-})\cap v=\emptyset $, and since $\omega
(B^{-}) $ is $\widetilde{f}$-invariant, 
$$
\widetilde{f}(\Theta _1)\cap v=\Theta _1\cap \widetilde{f}(v)=\emptyset . 
$$
Also, by the choice of $\epsilon >0$$,\widetilde{f}(v)\cap v=\emptyset .$

\begin{proposition}
\label{fimmesmo}: The following inclusion holds: $\Omega \subset \tilde
f(\Omega )$
\end{proposition}

{\it Proof:\ }

There are 2 possibilities:

\begin{enumerate}
\item  {$\tilde f(\Theta _1)\neq \Theta _1$} $\Rightarrow $ {$\tilde
f(\Theta _1)\cap \Theta _1=\emptyset $}

\item  {$\tilde f(\Theta _1)=\Theta _1$}
\end{enumerate}

Assume first that $\tilde f(\Theta _1)\cap \Theta _1=\emptyset $. Then 
$$
\partial \widetilde{f}(\Omega )\cap \partial \Omega =\emptyset . 
$$
Since $]-\infty ,a[\times \{0\}\subset \Omega $ and $\widetilde{f}(]-\infty
,a[\times \{0\})\supset ]-\infty ,a[\times \{0\},\Omega \cap \widetilde{f}%
(\Omega )\neq \emptyset .$ As $\left( \Theta _1\cup v\right) \cap \widetilde{%
f}(v)=\emptyset $ and $\widetilde{f}(v)\cap \overline{\Omega }^c\neq
\emptyset ,$ we get that $\widetilde{f}(v)\cap \overline{\Omega }=\emptyset
. $ And this implies that $\widetilde{f}(\Theta _1)\cap \overline{\Omega }%
=\emptyset ,$ because we are assuming that $\tilde f(\Theta _1)\cap \Theta
_1=\emptyset .$ So, if $\widetilde{z}\in \Omega ,$ there is a simple
continuous arc $\alpha \subset \Omega $ which connects $\widetilde{z}$ to
some point $\widetilde{z}_0\in ]-\infty ,a[\times \{0\}.$ As $\widetilde{z}%
_0\in \widetilde{f}(\Omega )$ and $\alpha \cap \left( \widetilde{f}(\Theta
_1)\cup \widetilde{f}(v)\right) =\emptyset ,$ we get that $\alpha \subset 
\widetilde{f}(\Omega )$ and so $\widetilde{f}(\Omega )\supset \Omega .$

Now, suppose that $\tilde f(\Theta _1)=\Theta _1.$ This implies that%
$$
\tilde f(\Omega )\text{ is a connected component of }\left( \Theta _1\cup
\tilde f(v)\right) ^c, 
$$
and 
$$
]-\infty ,a[\times \{0\}\subset \widetilde{f}(]-\infty ,a[\times
\{0\})\subset \tilde f(\Omega ). 
$$

As $\left( \Theta _1\cup v\right) \cap \tilde f(v)=\emptyset ,$ $\tilde f(v)$
does not intersect $\partial \Omega .$ Since both $v$ and $\Omega $ are
connected and $\tilde f(a,0)\in \widetilde{f}(v)$ does not belong to $\Omega
,$ we get that $\tilde f(v)\cap \Omega =\emptyset .$

Now, as above, let $\widetilde{z}$ be a point in $\Omega $ and $\alpha $ be
a simple continuous arc contained in $\Omega $ connecting $\widetilde{z}$ to
some $\widetilde{z}_0\in ]-\infty ,a[\times \{0\}.$ Since $\alpha \cap
\Theta _1=\alpha \cap \tilde f(v)=\emptyset ,$ $\alpha $ is contained in a
connected component of $(\Theta _1\cup \tilde f(v))^c.$ And since $
\widetilde{z}_0\in \alpha \cap \tilde f(\Omega ),$ it follows that $\alpha
\subset \tilde f(\Omega ).$ But this shows that any point $\widetilde{z}\in
\Omega $ is a point of $\tilde f(\Omega ),$ that is, $\Omega \subset \tilde
f(\Omega ).$ $\Box $

\vskip 0.2truecm

As $\Omega $ is open, the transitivity of $\widetilde{f}$ and the last
proposition yields that it is dense in the strip. But $\Omega \subset V^{-},$
arriving in the final contradiction that proves theorem 2. This same
argument is used often in the proofs of the next theorems.

\section{Proof of theorem 1}

Assume by contradiction that $\overline{p(B^{-})}\neq A.$ From fact \ref
{importante}, we know that there exists a connected component $\Gamma $ of $%
B^{-}$ that satisfies: $p(\Gamma )\subset \gamma _E^{-},$ $\overline{%
p(\Gamma )}\supset S^1\times \{0\},$ see expression (\ref{defgamaE}) and
figure 6. Let $\sigma >0$ be the number defined in the previous proof.
Let $\epsilon >0$ be sufficiently small such that:

\begin{itemize}
\item  $S^1\times [0,\epsilon ]\subset \gamma _E^{-};$

\item  for all $(\widetilde{x},\widetilde{y})\in {\rm I}\negthinspace {\rm R}%
\times [0,\epsilon ],$ $p_1\circ \tilde f(\tilde x,\widetilde{y})>\tilde
x+\sigma ;$
\end{itemize}

As $\overline{p(\Gamma )}\supset S^1\times \{0\},$ there exists a
sufficiently negative $a$ such that 
\begin{equation}
\label{definia2}\Gamma \cap \{a\}\times [0,\epsilon ]\neq \emptyset . 
\end{equation}

As in the previous proof, $B^{-}$ is closed, so there must be a $\delta \le
\epsilon $ such that $(a,\delta )\in B^{-}$, and such that for all $0\le
\tilde y<\delta ,(a,\tilde y)\notin B^{-},$ that is, $(a,\delta )$ is the
``lowest'' point of $B^{-}$ in $\{a\}\times [0,\epsilon ].$ We again denote
by $v$ the segment $\{a\}\times [0,\delta [.$

Let $\Gamma _1$ be the connected component of $B^{-}$ that contains $%
(a,\delta ).$ As $p(\Gamma _1)\cap \gamma _E^{-}\neq \emptyset ,$ we get
that $\overline{p(\Gamma _1)}\subset \gamma _E^{-},$ something that implies
the following important facts:

$$
\begin{array}{c}
dist(p(\Gamma _1),\ S^1\times \{1\})>0 \\ 
dist(\Gamma _1,{\rm I}\negthinspace {\rm R}\times \{1\})>0 
\end{array}
$$

We claim that $\Gamma _1{}^{+}$ is not above $\Gamma _1.$ We need the
following propositions:

\begin{proposition}
\label{primfabio}: $\Gamma _1{}^{+}\cap v=\emptyset .$
\end{proposition}

{\it Proof:}

This follows from $\Gamma^{+}\subset B^{-}$ since, by the definition of $v$, 
$B^{-}\cap v =\emptyset.$ $\Box $

\vskip 0.2truecm

\begin{proposition}
\label{segfabio}:\ If $\Gamma _1\cap \widetilde{f}(\Gamma _1)=\emptyset $
and $\Gamma _1\prec \widetilde{f}(\Gamma _1),$ then $\widetilde{f}(v)\cap
\Gamma _1=\emptyset .$
\end{proposition}

{\it Proof: }

As $\widetilde{f}(\Gamma _1)\subset B^{-},$ either $\Gamma _1\supset 
\widetilde{f}(\Gamma _1)$ or $\Gamma _1\cap \widetilde{f}(\Gamma
_1)=\emptyset .$ So, if $\Gamma _1\prec \widetilde{f}(\Gamma _1),$ we get by
lemma \ref{presordem} that $\widetilde{f}^{-1}(\Gamma _1)\prec \Gamma _1,$
so $\left[ \widetilde{f}^{-1}(\Gamma _1)\cap V_a^{-}\right] \subset \Gamma
_{1a,down}.$

On the other hand, note that $\Gamma _1\cup v$ is a closed connected set and 
$(\Gamma _1\cup v)^c$ has a connected component, denoted $\Omega ,$ which
contains $]-\infty ,a[\times \{0\}$ and another one which contains $%
]a,+\infty [\times \{0\}\cup {\rm I}\negthinspace {\rm R}\times \{1\}$.
Moreover, $\Omega \subset p^{-1}(\gamma _E^{-}).$ Also, it is immediate to
see that 
$$
closure(\Omega )\supset closure(\Gamma _{1a,down}), 
$$
so as $\widetilde{f}^{-1}(\Gamma _1)\cap \Gamma _1=\emptyset $ and $\left[ 
\widetilde{f}^{-1}(\Gamma _1)\cap V_a^{-}\right] \subset \Gamma
_{1a,down}\subset closure(\Omega ),$ we have two possibilities:

i) $\widetilde{f}^{-1}(\Gamma _1)\cap v=\emptyset ,$ something that implies
the proposition;

ii) $\widetilde{f}^{-1}(\Gamma _1)\cap v\neq \emptyset .$ 
Consider an element $\Theta \in \left[ \widetilde{f}^{-1}(\Gamma _1)\cap
closure(\Omega )\right] =\{$unlimited connected components of $\widetilde{f}%
^{-1}(\Gamma _1)\cap closure(\Omega )\}.$ The connectivity of $\widetilde{f}%
^{-1}(\Gamma _1)$ and the fact that $\widetilde{f}^{-1}(\Gamma _1)\cap
\Gamma _1=\emptyset $ imply that $\Theta $ intersects $v.$ As $\Theta
\subset V^{-}$ and $\widetilde{f}(\Theta )\subset \Gamma _1,$ we get that $%
\Theta \subset B^{-},$ something in contradiction with $B^{-}\cap
v=\emptyset .$ $\Box $

\vskip 0.2truecm

Moreover, if $\Gamma _1\prec \Gamma _1{}^{+},$ then $\left[ \widetilde{f}%
(\Gamma _1)\cap V_a^{-}\right] \subset $ $\left[ \Gamma _1^{+}\cap
V_a^{-}\right] \subset \Gamma _{1a,up},$ which implies, by the proof of
lemma \ref{partecha3}, that $\Gamma _1\prec \widetilde{f}(\Gamma _1).$

\begin{lemma}
\label{terfabio}:$\;\tilde f(\Gamma _1)$ is not above $\Gamma _1,$ that is
either $\Gamma _1\supset \widetilde{f}(\Gamma _1)$ or $\widetilde{f}(\Gamma
_1)\prec \Gamma _1,$ which implies that either $\Gamma _1=\Gamma _1{}^{+},$
or $\Gamma _1{}^{+}$ $\prec \Gamma _1{}.$
\end{lemma}

{\it Proof:}

Suppose $\Gamma _1\prec \tilde f(\Gamma _1).$ As in the previous
proposition, let $\Omega $ be the open connected component of $(\Gamma
_1\cup v)^c$ that is unlimited to the left and lies in $]-\infty ,0[\times
[0,1]\cap p^{-1}(\gamma _E^{-}).$ Clearly, $\partial \Omega \subset \Gamma
_1\cup v$ and as $\tilde f(\Gamma _1)\cap \Gamma _1=\emptyset ,$ $v\cap
\tilde f(v)=\emptyset $ and $\tilde f(\Gamma _1)\cap v=\tilde f(v)\cap
\Gamma _1=\emptyset ,$ we obtain 
$$
\tilde f(\partial \Omega )\cap \partial \Omega =\emptyset 
$$
because $\tilde f(\partial \Omega )=\partial \tilde f(\Omega )\subset $ $%
\tilde f(\Gamma _1)\cup \tilde f(v),$ which is a closed connected set that
does not intersect $(\Gamma _1\cup v)\supset \partial \Omega .$ As $\tilde
f(v)\cap \overline{\Omega }^c\neq \emptyset ,$ because $]-\infty ,a[\times
\{0\}\subset \tilde f(]-\infty ,a[\times \{0\}),$ we get that $\left( \tilde
f(\Gamma _1)\cup \tilde f(v)\right) \cap \overline{\Omega }=\emptyset $ and
so $\Omega \subset \tilde f(\Omega ).$ But, since $\Omega $ is open and
limited to the right, this contradicts the transitivity of $\widetilde{f}.$ $%
\Box $

\vskip 0.2truecm

So, either $\Gamma _1=\Gamma _1{}^{+},$ or $\Gamma _1{}^{+}$ $\prec \Gamma
_1{}.$ Two different cases may arise.

\subsection{$\Gamma _1$ is an injective component}

From lemma \ref{maisimp}, as $\overline{p(\Gamma _1)}\subset \gamma _E^{-},$
we obtain that $\overline{p(\Gamma _1)}\supset S^1\times \{0\}.$ So $%
dist(\Gamma _1,$ ${\rm I}\negthinspace {\rm R}\times \{0\})=0$ and lemma \ref
{acumula} implies that $\Gamma _1$ is a down component.

Lemma \ref{terfabio} gives two possibilities:

\begin{enumerate}
\item  $\tilde f(\Gamma _1)\subset \Gamma _1,$ that is, $\Gamma
_1{}^{+}=\Gamma _1$

\item  $\Gamma _1{}^{+}\prec \Gamma _1$
\end{enumerate}

\vskip 0.1truecm

If $\Gamma _1{}^{+}=\Gamma _1,$ then $\Gamma _1{}^{+}$ is also an injective
down component and for all integers $k>0,$ $\Gamma _1{}^{+}+(k,0)\prec
\Gamma _1{}.$

If $\Gamma _1{}^{+}\prec \Gamma _1,$ then we can prove the following:

\begin{fact}
\label{igantes}: For all integers $k>0,$ $\Gamma _1{}^{+}+(k,0)\cap \Gamma
_1=\emptyset $ and $\Gamma _1{}^{+}+(k,0)\prec \Gamma _1.$
\end{fact}

{\it Proof:}

If for some integer $k_0>0,$ $\Gamma _1{}^{+}+(k_0,0)$ intersects $\Gamma
_1, $ then ${}\Gamma _1{}^{+}$ intersects $\Gamma _1-(k_0,0).$ As $\Gamma
_1{}^{+}$ is a connected component of $B^{-}\;$and $\Gamma _1-(k_0,0)\subset
B^{-}$ is closed and connected, we get that $\Gamma _1{}^{+}\supset \Gamma
_1-(k_0,0).$ But this contradicts $\Gamma _1{}^{+}\prec \Gamma _1$ because
as $\Gamma _1$ is a down component, $\Gamma _1\prec \Gamma _1-(k,0).$ So for
all positive integers $k,$ $\Gamma _1{}^{+}+(k,0)$ does not intersect $%
\Gamma _1.$ If the fact is not true, then lemma \ref{partecha3} implies that
for some $k_{*}>0,$ $\Gamma _1{}\prec \Gamma _1{}^{+}+(k_{*},0),$ which
implies that $\Gamma _1-(k_{*},0)\prec \Gamma _1{}^{+}\prec \Gamma _1{}$ and
this again contradicts the fact that $\Gamma _1$ is a down component. $\Box $

\vskip 0.2truecm

So, in cases 1 and 2 above, for all integers $k>0,$ $\Gamma
_1{}^{+}+(k,0)\cap \Gamma _1=\emptyset $ and $\Gamma _1{}^{+}+(k,0)\prec
\Gamma _1.$

The important result of this subsection is the following:

\begin{lemma}
\label{omegalim}: There exists a vertical $V_r=\{r\}\times [0,1]$\ and a
sequence $n_i\stackrel{i\rightarrow \infty }{\rightarrow }\infty $ such that 
$\widetilde{f}^{n_i}(\Gamma _1)\cap V_r\neq 0$ for all $i.$
\end{lemma}

{\it Proof:}

As $\Gamma _1$ is a down component, $\Gamma _1\prec \Gamma _1-(1,0).$ Lemma 
\ref{presordem} tell us that $\widetilde{f}(\Gamma _1)\prec \widetilde{f}%
(\Gamma _1)-(1,0)$. Note that $\widetilde{f}(\Gamma _1)$ may not be a whole
connected component of $B^{-},$ but we will abuse notation and say that $
\widetilde{f}(\Gamma _1)$ is a down component.

Above we proved that $\widetilde{f}(\Gamma _1)+(k,0)\prec \Gamma _1$ for all
integers $k>0,$ so as $\widetilde{f}(\Gamma _1)\subset B^{-},$ in any of the
possibilities 1) or 2), $\widetilde{f}(\Gamma _1)\subset \Gamma _1\cup
\Omega $ because either:

\begin{enumerate}
\item  $\widetilde{f}(\Gamma _1)\subset \Gamma _1$

\item  $\widetilde{f}(\Gamma _1)\prec \Gamma _1$ $\Leftrightarrow $ $\left[ 
\widetilde{f}(\Gamma _1)\cap V_a^{-}\right] \subset closure(\Gamma
_{1a,down})$ (see expression (\ref{definia2}) for a definition of $a$). As $%
closure(\Gamma _{1a,down})\subset closure(\Omega ),$ which is a connected
set (see the proof of proposition \ref{segfabio} for a definition of $\Omega 
$) and $\partial \Omega \subset \Gamma _1\cup v$ does not intersect $
\widetilde{f}(\Gamma _1),$ we get that $\widetilde{f}(\Gamma _1)\subset
\Omega .$
\end{enumerate}

Let us fix some $k^{\prime }>0$ in a way that $\widetilde{f}(\Gamma
_1)+(k^{\prime },0)$ intersects $v=\{a\}\times [0,\delta [,$ see (\ref
{definia2})$.$ The reason why such a $k^{\prime }$ exists is the following:
As $\widetilde{f}(\Gamma _1)+(k,0)\cap \Gamma _1=\emptyset $ and $\widetilde{%
f}(\Gamma _1)+(k,0)\prec \Gamma _1$ for all integers $k>0,$ we get that $%
\left[ \widetilde{f}(\Gamma _1)+(k,0)\cap V_a^{-}\right] \subset
closure(\Gamma _{1a,down})\subset closure(\Omega ).$ And as $\partial \Omega
\subset \Gamma _1\cup v$ and $\overline{\Omega }\subset ]-\infty ,0]\times
[0,1],$ we get that if $k^{\prime }>0$ is sufficiently large in a way that $
\widetilde{f}(\Gamma _1)+(k^{\prime },0)$ intersects $\{1\}\times [0,1],$
then $\widetilde{f}(\Gamma _1)+(k^{\prime },0)$ intersects the boundary of $%
\Omega $ in the only possible place, $v.$ Denote by $\Gamma ^{*}$ an
unlimited connected component of $\widetilde{f}(\Gamma _1)+(k^{\prime
},0)\cap closure(\Omega ).$ By the choice of $k^{\prime }>0$ and the
connectivity of $\widetilde{f}(\Gamma _1)+(k^{\prime },0),$ we get that $%
\Gamma ^{*}$ is not contained in $B^{-}$ because it intersects $v.$ So,
there exists a positive integer $a_1>0$ such that $\widetilde{f}%
^{a_1}(\Gamma ^{*})$ intersects $]0,+\infty [\times [0,1].$ Remember that $
\widetilde{f}(\Gamma _1)+(k^{\prime },0)\prec \Gamma _1$ and, as we said
above, either $\widetilde{f}(\Gamma _1)\subset \Gamma _1$ or $\widetilde{f}%
(\Gamma _1)\prec \Gamma _1.$ In case $\widetilde{f}(\Gamma _1)\subset \Gamma
_1,$ we get that $\widetilde{f}^{a_1+1}(\Gamma _1)\subset \Gamma _1.$ Before
continuing the proof, let us state the following:

\begin{proposition}
\label{faltando}: If $\Gamma $ is a connected component of $B^{-}$ which
satisfies $\widetilde{f}(\Gamma )\cap \Gamma =\emptyset $ and $\widetilde{f}%
(\Gamma )\prec \Gamma ,$ then $\widetilde{f}^n(\Gamma )\cap \Gamma
=\emptyset $ and $\widetilde{f}^n(\Gamma )\prec \Gamma $ for all integers $%
n>0.$
\end{proposition}

{\it Proof:}

By contradiction, suppose there exists some $n_0>1$ (then smallest one) such
that $\widetilde{f}^{n_0}(\Gamma )\cap \Gamma \neq \emptyset .$ This means
that $\Gamma ,$ $\widetilde{f}(\Gamma ),\widetilde{f}^2(\Gamma ),...,$ $
\widetilde{f}^{n_0-1}(\Gamma )$ are disjoint closed connected subsets of the
strip $\widetilde{A},$ each of them having a connected complement and $
\widetilde{f}^{n_0}(\Gamma )\subset \Gamma $. As $\widetilde{f}(\Gamma
_1)\prec \Gamma _1,$ lemmas \ref{associati} and \ref{presordem} imply that 
\begin{equation}
\label{falt1}\widetilde{f}^{n_0-1}(\Gamma )\prec ...\prec \widetilde{f}%
^2(\Gamma )\prec \widetilde{f}(\Gamma )\prec \Gamma . 
\end{equation}
On the other hand, as $\widetilde{f}^{n_0}(\Gamma )\cap \widetilde{f}%
^{n_0-1}(\Gamma )=\emptyset ,$ lemma \ref{presordem} implies that $
\widetilde{f}^{n_0}(\Gamma )\prec \widetilde{f}^{n_0-1}(\Gamma ).$ So, $
\widetilde{f}^{n_0}(\Gamma )\cap \left( \widetilde{f}^{n_0-1}(\Gamma
)\right) _{down}$ has an unlimited connected component. As $\widetilde{f}%
^{n_0}(\Gamma )\subset \Gamma $ and $\widetilde{f}^{n_0-1}(\Gamma )\cap
\Gamma =\emptyset ,$ using lemma \ref{partecha2} we get that $\Gamma \prec 
\widetilde{f}^{n_0-1}(\Gamma ),$ a contradiction with expression (\ref{falt1}%
). So, $\widetilde{f}^n(\Gamma )\cap \Gamma =\emptyset $ for all integers $%
n>0.$ The other implication follows from lemma \ref{presordem}. $\Box $

\vskip 0.2truecm

So, if $\widetilde{f}(\Gamma _1)\prec \Gamma _1$ ($\widetilde{f}(\Gamma
_1)\cap \Gamma _1=\emptyset $), then $\widetilde{f}^{a_1+1}(\Gamma _1)\cap
\Gamma _1=\emptyset $ and $\widetilde{f}^{a_1+1}(\Gamma _1)\prec \Gamma _1.$
As $\Gamma _1$ is a down component, $\widetilde{f}^{a_1}(\widetilde{f}%
(\Gamma _1)+(k^{\prime },0))=\widetilde{f}^{a_1+1}(\Gamma _1)+(k^{\prime
},0)\prec \widetilde{f}^{a_1+1}(\Gamma _1).$ Now, note that $\widetilde{f}%
^{a_1+1}(\Gamma _1)+(k^{\prime },0)\cap \Gamma _1=\emptyset .$

This happens because, if $\widetilde{f}^{a_1+1}(\Gamma _1)+(k^{\prime
},0)\cap \Gamma _1\neq \emptyset ,$ then there would be a connected
component of $B^{-},$ denoted $\Psi ,$ containing both $\widetilde{f}%
^{a_1+1}(\Gamma _1)$ and $\Gamma _1-(k^{\prime },0).$ Clearly, $\Psi $ is
not $\Gamma _1$ and both $\Psi \cap \Gamma _{1down}$ and $\Psi \cap \Gamma
_{1up}$ have unlimited connected components, because $\widetilde{f}%
^{a_1+1}(\Gamma _1)\prec \Gamma _1$ and $\Gamma _1\prec \Gamma _1-(k^{\prime
},0),$ a contradiction. Thus, $\widetilde{f}^{a_1}(\widetilde{f}(\Gamma
_1)+(k^{\prime },0))=\widetilde{f}^{a_1+1}(\Gamma _1)+(k^{\prime },0)\prec
\Gamma _1$ and so, as

$$
\Gamma ^{*}\subset \widetilde{f}(\Gamma _1)+(k^{\prime },0)\ \Rightarrow 
\text{ }\widetilde{f}^{a_1}(\Gamma ^{*})\cap \Gamma _1=\emptyset \text{ and }
\widetilde{f}^{a_1}(\Gamma ^{*})\prec \Gamma _1. 
$$
As above, the fact that $\widetilde{f}^{a_1}(\Gamma ^{*})$ intersects $%
]0,+\infty [\times [0,1]$ implies that $\widetilde{f}^{a_1}(\Gamma ^{*})\cap
closure(\Omega )$ has an unlimited connected component, $\Gamma ^{**}$ which
intersects $v.$ So, $\Gamma ^{**}$ is not contained in $B^{-}$ and thus
there exists an integer $a_2>0$ such that $\widetilde{f}^{a_2}(\Gamma
^{**})\subset \widetilde{f}^{a_2+a_1+1}(\Gamma _1)+(k^{\prime },0)$
intersects $]0,+\infty [\times [0,1].$ In exactly the same way as above, we
obtain an unlimited connected component of $\widetilde{f}^{a_2}(\Gamma
^{**})\cap closure(\Omega )$, denoted $\Gamma ^{***}$ which intersects $v.$
So, $\Gamma ^{***}$ is not contained in $B^{-}$ and there exists an integer $%
a_3>0$ such that $\widetilde{f}^{a_3}(\Gamma ^{***})$ intersects $]0,+\infty
[\times [0,1]$ and so on.

Thus, if we define $n_i=a_1+a_2+...+a_i+1,$ we get that $n_i\stackrel{%
i\rightarrow \infty }{\rightarrow }\infty $ and for all $i\geq 1,$ $
\widetilde{f}^{n_i-1}(\widetilde{f}(\Gamma _1)+(k^{\prime },0))\supset 
\widetilde{f}^{n_i-1}(\Gamma ^{*})=\widetilde{f}^{a_2+...+a_i}(\widetilde{f}%
^{a_1}(\Gamma ^{*}))\supset \widetilde{f}^{a_2+...+a_i}(\Gamma ^{**})\supset
...\supset \widetilde{f}^{a_i}(\Gamma ^{\stackrel{i-times}{*...*}})$ and $
\widetilde{f}^{a_i}(\Gamma ^{\stackrel{i-times}{*...*}})$ intersects $%
]0,+\infty [\times [0,1].$ So, 
$$
\widetilde{f}^{n_i}(\Gamma _1)\text{ intersects }V_0-(k^{\prime
},0)=V_{-k^{\prime }} 
$$
and the lemma is proved. $\Box $

\vskip 0.2truecm

But the lemma implies that $\omega(B^{-})\neq\emptyset,$ contradicting
theorem 2. Therefore we must have that $\Gamma_1$ is a non-injective
component.

\subsection{$\Gamma _1$ is an non-injective component}

From lemma \ref{maisimp}, as $\overline{p(\Gamma _1)}\subset \gamma _E^{-},$
we obtain that $\overline{p(\Gamma _1)}\supset S^1\times \{0\}.$

The next result is interesting by itself:

\begin{fact}
\label{fmenos}:\ If $\Gamma $ is a non-injective component of $B^{-}$ and $
\widetilde{f}(\Gamma )\subset \Gamma $, then there exists an integer $k>0$
such that $\widetilde{f}^{-1}(\Gamma )\subset \Gamma +(k,0).$
\end{fact}

{\it Proof: }

Since $\tilde f(\Gamma )\subset \Gamma ,$ we get $\Gamma \subset \tilde
f^{-1}(\Gamma ).$ As $\Gamma \subset V^{-},$ $\tilde f^{-1}(\Gamma )$ is
limited to the right, so there exists an integer $k>0$ such that $\tilde
f^{-1}(\Gamma )-(k,0)\subset V^{-}.$ If $i\geq 1,$%
$$
\widetilde{f}^i\left( \tilde f^{-1}(\Gamma )-(k,0)\right) =\tilde
f^{i-1}(\Gamma )-(k,0)\subset \Gamma -(k,0)\subset \Gamma . 
$$
So, as the closed connected set $\tilde f^{-1}(\Gamma )-(k,0)$ is unlimited
to the left and has all its positive iterates in $V^{-}$, it is contained in 
$B^{-}.$ As $\Gamma \subset \tilde f^{-1}(\Gamma )$ $\Rightarrow $ $\Gamma
-(k,0)\subset \tilde f^{-1}(\Gamma )-(k,0).$ So, $\tilde f^{-1}(\Gamma
)-(k,0)$ intersects $\Gamma $ because $\Gamma \supset \Gamma -(k,0).$ But $%
\Gamma $ is a connected component of $B^{-}$ and $\widetilde{f}^{-1}(\Gamma
)-(k,0)$ is connected, therefore 
$$
\tilde f^{-1}(\Gamma )-(k,0)\subset \Gamma , 
$$
something that proves the fact. $\Box $

\vskip 0.2truecm

Clearly, for any integer $n\geq 1,$%
$$
\widetilde{f}^{-n}(\Gamma )\subset \Gamma +(n.k,0). 
$$

In contrast with the case when $\Gamma _1$ is injective, lemma \ref{terfabio}
implies that the only possibility here is $\tilde f(\Gamma _1)\prec \Gamma
_1 $ because of the next lemma:

\begin{lemma}
\label{noninj1}: It is not possible that $\tilde f(\Gamma _1)\subset \Gamma
_1.$
\end{lemma}

{\it Proof: }

Suppose that $\tilde f(\Gamma _1)\subset \Gamma _1.$ Let 
$$
\Omega _1=\stackrel{\infty }{\stackunder{n=0}{\cup }}\widetilde{f}%
^{-n}(\Omega ), 
$$
where, as in the proofs of proposition \ref{segfabio} and lemma \ref
{terfabio}, $\Omega $ is the open connected component of $(\Gamma _1\cup
v)^c $ that contains $]-\infty ,a[\times\{0\}.$ Note that $\Omega$ is 
contained in $p^{-1}(\gamma _E^{-})$ and clearly, $\widetilde{f}^{-1}(\Omega
_1)\subset \Omega _1.$ Moreover, the following is true:

\begin{proposition}
\label{contido}: $\Omega _1$ is contained in $p^{-1}(\gamma _E^{-}).$
\end{proposition}

Before proving this proposition, let us show how it is used to prove our
lemma. Since $\Omega_1$ is open and $\widetilde{f}^{-1}(\Omega_1)\subset
\Omega_1$, we must have, by the transitivity of $\widetilde{f},$ that $%
\Omega_1$ is dense. But this contradicts the proposition. $\Box $

\vskip 0.2truecm

{\it Proof of proposition }\ref{contido}:

First note that, as the boundary of $\Omega $ is contained in $\Gamma _1\cup
v$, for all integers $i>0$ we have: 
\begin{equation}
\label{sobrbordo}\partial \left( \widetilde{f}^{-i}(\Omega )\right) \subset 
\stackrel{\infty }{\stackunder{n=0}{\cup }}\widetilde{f}^{-n}(\Gamma _1\cup
v)=\left( \stackrel{\infty }{\stackunder{n=0}{\cup }}\widetilde{f}%
^{-n}(\Gamma _1)\right) \cup \left( \stackrel{\infty }{\stackunder{n=0}{\cup 
}}\widetilde{f}^{-n}(v)\right) . 
\end{equation}

Clearly $\Omega _1$ is an open set. Let us show that it is connected. Each
set of the form $\widetilde{f}^{-i}(\Omega )$ is connected because $
\widetilde{f}$ is a homeomorphism. Also, since $\widetilde{f}^{-i}(]-\infty
,a[\times \{0\})\subset ]-\infty ,a[\times \{0\}$, we have $\widetilde{f}%
^{-i}(\Omega )\cap \Omega \not =\emptyset .$ But $\Omega $ is also open and
connected, so $\Omega _1$ must be connected.

For all integers $i>0,$ as $\widetilde{f}^{-i}(\Omega )$ is connected,
intersects ${\rm I}\negthinspace {\rm R}\times \{0\}$ and is disjoint from $%
{\rm I}\negthinspace {\rm R}\times \{1\},$ if we show that $\left( \stackrel{%
\infty }{\stackunder{n=0}{\cup }}\widetilde{f}^{-n}(\Gamma _1)\right) \cup
\left( \stackrel{\infty }{\stackunder{n=0}{\cup }}\widetilde{f}%
^{-n}(v)\right) \subset p^{-1}(\gamma _E^{-}),$ then expression (\ref
{sobrbordo}) implies that $\widetilde{f}^{-i}(\Omega )\subset p^{-1}(\gamma
_E^{-}),$ which gives:\ $\Omega _1\subset p^{-1}(\gamma _E^{-})$ and the
proof is complete.

Let us analyze first what happens to $\widetilde{f}^{-n}(\Gamma _1),$ for
all integers $n>0.$

From fact \ref{fmenos}, $\stackrel{\infty }{\stackunder{n=0}{\cup }} 
\widetilde{f}^{-n}(\Gamma _1)\subset \stackrel{\infty }{\stackunder{n=0}{%
\cup }}(\Gamma _1+(n,0))\subset p^{-1}(\gamma _E^{-}).$

We are left to deal with $\stackrel{\infty }{\stackunder{n=0}{\cup }} 
\widetilde{f}^{-n}(v).$ Let us show that $\widetilde{f}^{-1}(v)\subset
\Omega .$ From the choice of $v,$ $\widetilde{f}^{-1}(v)\cap v=\emptyset .$
Also, from the definition of $v=\{a\}\times [0,\delta [,$ we get that $
\widetilde{f}^{-1}(v)\cap \Gamma _1=\emptyset $ because $v\cap
B^{-}=\emptyset .$ Finally, the following inclusions 
$$
\Omega \supset ]-\infty ,a[\times \{0\}\text{ and }]-\infty ,a[\times
\{0\}\supset \widetilde{f}^{-1}(]-\infty ,a[\times \{0\})\text{ } 
$$
imply that $\widetilde{f}^{-1}(v)\cap \Omega \neq \emptyset $ $\Rightarrow $ 
$\widetilde{f}^{-1}(v)\subset \Omega \subset p^{-1}(\gamma _E^{-}).$

So, $\widetilde{f}^{-2}(v)\subset \widetilde{f}^{-1}(\Omega ),$ whose
boundary, $\partial \left( \widetilde{f}^{-1}(\Omega )\right) ,$ is
contained in

$\left( \stackrel{\infty }{\stackunder{n=0}{\cup }}(\Gamma _1+(n,0))\right)
\cup \widetilde{f}^{-1}(v)\subset p^{-1}(\gamma _E^{-}).$ As above, as $
\widetilde{f}^{-1}(\Omega )$ is connected, intersects ${\rm I}%
\negthinspace 
{\rm R}\times \{0\}$ and is disjoint from ${\rm I}\negthinspace {\rm R}%
\times \{1\},$ we get that 
$$
\widetilde{f}^{-2}(v)\subset \widetilde{f}^{-1}(\Omega )\subset
p^{-1}(\gamma _E^{-}). 
$$
So, $\widetilde{f}^{-3}(v)\subset \widetilde{f}^{-2}(\Omega )$ and an
analogous argument implies that $\widetilde{f}^{-3}(v)\subset \widetilde{f}%
^{-2}(\Omega )\subset p^{-1}(\gamma _E^{-}).$ An induction shows that%
$$
\widetilde{f}^{-n}(v)\subset \widetilde{f}^{-n+1}(\Omega )\subset
p^{-1}(\gamma _E^{-})\text{ for all integers }n\geq 1, 
$$
and the proposition is proved. $\Box $

\vskip 0.2truecm

Thus, if $\Gamma_1$ is a non-injective component, the only possibility is $%
\tilde f(\Gamma _1)\prec \Gamma _1.$ The next lemma is a version of lemma 
\ref{omegalim} to a non-injective $\Gamma _1:$

\begin{lemma}
\label{omegalim2}: There exists a vertical $V_r=\{r\}\times [0,1]$\ and a
sequence $n_i\stackrel{i\rightarrow \infty }{\rightarrow }\infty $ such that 
$\widetilde{f}^{n_i}(\Gamma _1)\cap V_r\neq 0$ for all integer $i\geq 1.$
\end{lemma}

{\it Proof:}

As $\widetilde{f}(\Gamma _1)\prec \Gamma _1$ $\Rightarrow $ $\left[ 
\widetilde{f}(\Gamma _1)\cap V_a^{-}\right] \subset closure(\Gamma
_{1a,down})$ (see expression (\ref{definia2}) for a definition of $a$). As $%
closure(\Gamma _{1a,down})\subset closure(\Omega ),$ which is a connected
set and $\partial \Omega \subset \Gamma _1\cup v$ does not intersect $
\widetilde{f}(\Gamma _1),$ we get that $\widetilde{f}(\Gamma _1)\subset
\Omega .$

If for some integers $n_0>0$ and $k_0>0,$ $\widetilde{f}^{n_0}(\Gamma
_1)+(k_0,0)\cap \Gamma _1\neq \emptyset ,$ then $\widetilde{f}^{n_0}(\Gamma
_1)\cap \Gamma _1-(k_0,0)\neq \emptyset $ $\Rightarrow $ $\widetilde{f}%
^{n_0}(\Gamma _1)\cap \Gamma _1\neq \emptyset $ $\Rightarrow $ $\widetilde{f}%
^{n_0}(\Gamma _1)\subset \Gamma _1$, which, using proposition \ref{faltando}%
, implies that $\widetilde{f}(\Gamma _1)\subset \Gamma _1,$ a contradiction.
So, for all integers $n>0$ and $k>0,$ as $\widetilde{f}^n(\Gamma
_1)+(k,0)\supset \widetilde{f}^n(\Gamma _1)$ and 
$$
\tilde f^n(\Gamma _1)\prec ...\prec \tilde f(\Gamma _1)\prec \Gamma _1\text{
(see proposition \ref{faltando}),} 
$$
we get that $\widetilde{f}^n(\Gamma _1)+(k,0)\cap \Gamma _1=\emptyset $ and $
\widetilde{f}^n(\Gamma _1)+(k,0)\prec \Gamma _1.$

Now the proof goes exactly as in lemma \ref{omegalim}. $\Box $

\vskip 0.2truecm

Of course, we have arrived at the same contradiction as in subsection 4.1,
and so theorem 1 is proved.

\section{Proof of theorem 3}

\begin{lemma}
:\ There exists an integer $N_1>0$ such that $\tilde f^{N_1}(B^{-})\subset
B^{-}-(1,0)$
\end{lemma}

{\it Proof:}

Theorem 2 shows that $\omega (B^{-})=\emptyset ,$ so there must be an
integer $N_1>0$ such that, for all $n\ge N_1,\tilde f^n(B^{-})\subset
]-\infty ,-1[\times [0,1].$ Suppose, by contradiction, that there exists $
\widetilde{z}\in B^{-}$ such that $\tilde f^{N_1}(\widetilde{z})+(1,0)\notin
B^{-}.$

Let $\Gamma $ be the connected component of $B^{-}$ that contains $
\widetilde{z}.$ Clearly $\tilde f^{N_1}(\Gamma )\subset V_{-1}^{-}.$ As $%
f^{N_1}(\widetilde{z})+(1,0)\in \tilde f^{N_1}(\Gamma )+(1,0),$ $\tilde
f^{N_1}(\Gamma )+(1,0)$ is not a subset of $B^{-}.$ But $\tilde
f^{N_1}(\Gamma )+(1,0)$ is connected, unbounded and $\tilde f^{N_1}(\Gamma
)+(1,0)\subset V^{-},$ therefore there must be a $N_2>0$ such that $\tilde
f^{N_2}(\tilde f^{N_1}(\Gamma )+(1,0))=\tilde f^{N_1+N_2}(\Gamma )+(1,0)$ is
not contained in $V^{-}.$ This implies that $f^{N_1+N_2}(\Gamma )$ is not
contained in $V_{-1}^{-},$ a contradiction that proves the lemma. $\Box $

\vskip 0.2truecm

As $\tilde f^{N_1}(B^{-})\subset B^{-}-(1,0)$, for any positive integer $k,$ 
$$
\tilde f^{kN_1}(B^{-})\subset B^{-}-(k,0)\subset V_{-k}^{-}, 
$$
and so it follows that, for any point $\widetilde{z}\in B^{-},$ 
$$
\stackunder{n\rightarrow \infty }{\lim \sup }\frac{p_1(\tilde f^n(\widetilde{%
z}))-p_1(\widetilde{z})}n\le -\frac 1{N_1}, 
$$
and this proves theorem 3.

\section{Proof of theorem 4}

Let $\epsilon >0$ be such that for all $(\widetilde{x},\widetilde{y})\in 
{\rm I}\negthinspace {\rm R}\times \left\{ [0,\epsilon ]\cup [1-\epsilon
,1]\right\} ,$ $p_1\circ \tilde f(\tilde x,\widetilde{y})>\tilde x+\sigma ,$
for a certain fixed $\sigma >0.$ %










As theorem 1 says that $\overline{p(B^{-})}=A,$ there exists a sufficiently
negative $b$ such that 
$$
\Theta \cap \{b\}\times [0,\epsilon ]\neq \emptyset , 
$$
for some connected component $\Theta $ of $B^{-}.$ As in the beginning of
the proof of theorem 1, in the following we will consider the ``lowest''
component of $B^{-}$ in $\{b\}\times [0,\epsilon ].$

First, remember that as $B^{-}$ is closed, there must be a $0<\delta \le
\epsilon $ such that $(b,\delta )\in B^{-}$, and for all $0\le \tilde
y<\delta ,(b,\tilde y)\notin B^{-},$ that is, $(b,\delta )$ is the
``lowest'' point of $B^{-}$ in $\{b\}\times [0,\epsilon ].$ We denote by $v$
the segment $\{b\}\times [0,\delta [.$

Let $\Gamma _1$ be the connected component of $B^{-}$ that contains $%
(b,\delta ).$ By propositions \ref{primfabio}, \ref{segfabio} and lemma \ref
{terfabio}, if $\Gamma _1\prec \tilde f(\Gamma _1),$ then the set $\Omega ,$
which is the open connected component of $(\Gamma _1\cup v)^c$ that is
unlimited to the left, lies in $]-\infty ,0[\times [0,1]$ and contains $%
]-\infty ,b[\times \{0\}$ satisfies the following: $\Omega \subset \tilde
f(\Omega )$ $\Leftrightarrow $ $\widetilde{f}^{-1}(\Omega )\subset \Omega $
and this contradicts the existence of a dense orbit for $\widetilde{f}.$

So, either $\tilde f(\Gamma _1)\subset \Gamma _1$ or $\tilde f(\Gamma
_1)\prec \Gamma _1.$ In order to analyze the two previous possibilities, we
have to consider all possible ``shapes'' for $\Gamma _1:$

1) $\Gamma _1$ is an injective down component;

2) $\Gamma _1$ is an injective up component;

3)$\;\Gamma _1$ is a non-injective component;

In case 1, if $\tilde f(\Gamma _1)\subset \Gamma _1$ or $\tilde f(\Gamma
_1)\prec \Gamma _1$ and in case 3, if $\tilde f(\Gamma _1)\prec \Gamma _1$,
lemmas \ref{omegalim} and \ref{omegalim2} imply that $\omega (B^{-})$ (see (%
\ref{defomelim})) is not empty. And this is a contradiction with theorem 2.

So, either $\Gamma _1$ is an injective up component or $\Gamma _1$ is an
non-injective component and $\tilde f(\Gamma _1)\subset \Gamma _1.$ Suppose
that $\Gamma _1$ is an injective up component. We have two possibilities:

I) $dist(\Gamma _1,{\rm I}\negthinspace {\rm R}\times \{1\})>0;$

II)\ $dist(\Gamma _1,{\rm I}\negthinspace {\rm R}\times \{1\})=0;$

\begin{lemma}
\label{Inaoocor}: If $\Gamma _1$ is an injective up component, then $%
dist(\Gamma _1,{\rm I}\negthinspace {\rm R}\times \{1\})=0.$
\end{lemma}

{\it Proof:}

As $\Gamma _1$ is an injective up component, lemma \ref{acumula} implies
that 
$$
dist(\Gamma _1,{\rm I}\negthinspace {\rm R}\times \{0\})>0. 
$$
So if I) holds, there exists $\epsilon _1>0$ such that $\Gamma _1\cap {\rm I}%
\negthinspace {\rm R}\times \left\{ [0,\epsilon _1]\cup [1-\epsilon
_1,1]\right\} =\emptyset .$

Since $\tilde f$ is transitive, $f$ is transitive and thus there is a point $%
z\in S^1\times [1-\epsilon _1/2,1]$ and an integer $n>0$ such that $%
f^{-n}(z)\in S^1\times [0,\epsilon _1/2].$ We know that $\tilde f(\Gamma
_1)\subset \Gamma _1$ or $\tilde f(\Gamma _1)\prec \Gamma _1,$ so by
proposition \ref{faltando} and lemmas \ref{associati} and \ref{presordem} we
get that

$$
\tilde f^n(\Gamma _1)\subset \Gamma _1\ or\text{ }\ \tilde f^n(\Gamma
_1)\prec \Gamma _1. 
$$

Now let $d\in {\rm I}\negthinspace {\rm R}$ be such that $\widetilde{f}%
^{-i}(V_d)\subset V_{m_{\Gamma _1}-1}^{-}$ (see expression (\ref{mgama}))
for $i=0,1,...,n$ and $\tilde f^n(\Gamma _1)\cap V_d^{-}\subset \Gamma
_1\cup \Gamma _{1\text{ }down},$ see proposition \ref{limitado}.

Let $\widetilde{z}\in V_d^{-}$ be a point such that $p(\widetilde{z})=z$ and
let $k$ be the vertical line segment that has as extremes $\widetilde{z}$
and a point $\widetilde{z}_1$ in ${\rm I}\negthinspace {\rm R}\times \{1\}.$

As $\widetilde{f}^{-n}(\widetilde{z})\in {\rm I}\negthinspace {\rm R}\times
[0,\epsilon _1/2]\cap V_{m_{\Gamma _1}-1}^{-},$ we obtain that $\widetilde{f}%
^{-n}(\widetilde{z})\in \Gamma _{1\text{ }down}.$ As $\widetilde{f}%
^{-n}(V_d)\subset V_{m_{\Gamma _1}-1}^{-},$ we get that $\widetilde{f}%
^{-n}(k)\cap V_{m_{\Gamma _1}}=\emptyset .$ Since $\tilde f^{-n}(\widetilde{z%
}_1)\notin \Gamma _{1\text{ }down}$ and $\tilde f^{-n}(\tilde z)\in \Gamma
_{1\text{ }down},$ and since $k$ is connected, $\tilde f^{-n}(k)\cap
\partial (\Gamma _{1\text{ }down})\neq \emptyset .$ But $\partial (\Gamma
_{1 \text{ }down})\subset \Gamma _1\cup V_{m_{\Gamma _1}}$ and as $
\widetilde{f}^{-n}(k)\cap V_{m_{\Gamma _1}}=\emptyset ,$ we get that $
\widetilde{f}^{-n}(k)\cap \Gamma _1\neq \emptyset ,$ which implies that $%
k\cap \widetilde{f}^n(\Gamma _1)\neq \emptyset $ and this is a contradiction
because $k\subset V_d^{-}\cap \Gamma _{1\text{ }up}$ and $\tilde f^n(\Gamma
_1)\cap V_d^{-}\subset \Gamma _1\cup \Gamma _{1\text{ }down}.$ So I)\ does
not hold. $\Box $

\vskip 0.2truecm

Thus if $\Gamma _1$ is injective, then II)\ holds. So consider a
sufficiently negative $c$ such that 
$$
\Gamma _1\cap \{c\}\times [1-\epsilon ,1]\neq \emptyset , 
$$
where $\epsilon >0$ was defined in the beginning of this section. As above,
as $B^{-}$ is closed, there must be a $0<\mu \le \epsilon $ such that $%
(c,1-\mu )\in B^{-}$, and for all $1-\mu \le \tilde y<1,(c,\tilde y)\notin
B^{-},$ that is, $(c,1-\mu )$ is the ``highest'' point of $B^{-}$ in $%
\{c\}\times [1-\epsilon ,1].$ We denote by $w$ the segment $\{c\}\times
]1-\mu ,1].$

Let $\Gamma _2$ be the connected component of $B^{-}$ that contains $%
(c,1-\mu ).$ An argument analogous to the one which implies that $\Gamma _1$
can not be an injective down component, implies that $\Gamma _2$ can not be
an injective up component, so if $\Gamma _1$ is injective, $\Gamma _1\neq
\Gamma _2$ and thus $\Gamma _1\cap \Gamma _2=\emptyset .$ So $\Gamma _2$ is
either non-injective or an injective down component. In the second case, as $%
dist(\Gamma _2,{\rm I}\negthinspace {\rm R}\times \{1\})>0$ (see lemma \ref
{acumula}), it is not possible that $\Gamma _1\prec \Gamma _2.$ But this
implies that $\Gamma _2\prec \Gamma _1$ and so $\Gamma _2$ intersects $v.$
And this is a contradiction with the definition of $v.$ So $\Gamma _2$ is a
non-injective component. By exactly the same reasoning applied to $\Gamma
_1, $ we must have $\tilde f(\Gamma _2)\subset \Gamma _2.$

The following lemma concludes the proof of theorem 4, because either $\Gamma
_1$ is an non-injective component and $\tilde f(\Gamma _1)\subset \Gamma _1$
or, in case $\Gamma _1$ is an injective up component, $\Gamma _2$ is
non-injective and $\tilde f(\Gamma _2)\subset \Gamma _2$.

\begin{lemma}
\label{ramodenso}:\ If $\Gamma $ is a non-injective component of $B^{-}$
such that $\tilde f(\Gamma )\subset \Gamma ,$ then $\overline{p(\Gamma )}=A.$
\end{lemma}

{\it Proof:}

First of all, note that the set $\Gamma $ has all the properties required
for the set $D$ in subsection 2.2, so lemma \ref{conts1} implies that either 
$\overline{p(\Gamma )}\supset S^1\times \{0\}$ or $\overline{p(\Gamma )}%
\supset S^1\times \{1\}.$ So let us suppose, without loss of generality,
that 
$$
\overline{p(\Gamma )}\supset S^1\times \{0\}. 
$$

Lemma \ref{fullmeas} shows that, if $p(\Gamma )$ is not dense in $A,$ then
there exists a simple closed curve $\gamma \subset interior(A),$ which is
homotopically non trivial and such that $\overline{p(\Gamma )}\cap \gamma
=\emptyset .$ But since $p(\Gamma )$ is connected, we must have $\Gamma
\subset p^{-1}(\gamma ^{-}).$

As $\Gamma $ is closed and $S^1\times \{0\}\subset \overline{p(\Gamma )},$
we can find a point $(c^{\prime },\delta ^{\prime })\in \Gamma $ such that:

\begin{enumerate}
\item  $\delta ^{\prime }<\epsilon ,$ where $\epsilon >0$ was defined in the
beginning of this section;

\item  if $v^{\prime }=c^{\prime }\times [0,\delta ^{\prime }[,$ then $%
\Gamma \cap v^{\prime }=\emptyset $ and $\overline{v^{\prime }}\subset
p^{-1}(\gamma ^{-}).$
\end{enumerate}

Now, let us choose $\Omega ^{\prime }$ as the connected component of $\left(
\Gamma \cup v^{\prime }\right) ^c$ that contains $]-\infty ,c^{\prime
}[\times \{0\}$ and consider the following set, as we did in proposition \ref
{contido}:

$$
\Omega _{sat}=\stackrel{\infty }{\stackunder{n=0}{\cup }}\widetilde{f}%
^{-n}(\Omega ^{\prime }) 
$$

A simple repetition of the same arguments used in the proof of proposition 
\ref{contido} yields that $\Omega _{sat}\subset p^{-1}(\gamma ^{-}).$ Again,
since $\tilde f^{-1}(\Omega _{sat})\subset \Omega _{sat}$ this contradicts
the transitivity of $\tilde f$ and finishes the proof. $\Box $

\vskip 0.2truecm

As we know that at least one member of the set $\{\Gamma _1,\Gamma _2\}$ is
non-injective and positively invariant, the above lemma implies that $\Gamma
_1$ or $\Gamma _2$ must have a dense projection to the annulus. One more
thing can be said, which will be important in the proof of the next theorem:

\begin{proposition}
\label{todosmaisbai}: Both $\Gamma _1$ and $\Gamma _2$ are non-injective
components.
\end{proposition}

{\it Proof:}

\ If the proposition is not true for $\Gamma _1$, then as we already proved, 
$\Gamma _1$ must be an injective up component. As $\Gamma _2$ does not cross 
$v$ and $\Gamma _1$ does not cross $w,$ by the definitions of $v$ and $w,$
it must be the case that $\Gamma _1\prec \Gamma _2$ and so%
$$
dist(\Gamma _2,{\rm I}\negthinspace {\rm R}\times \{0\})>0, 
$$
something that contradicts lemma \ref{ramodenso}. So $\Gamma _1$ is a
non-injective component. To conclude the proof we have to note that $\Gamma
_1$ and $\Gamma _2$ have analogous properties, $\Gamma _1$ is the connected
component of $B^{-}$ that contains the lowest point of $B^{-}$ in $%
\{b\}\times [0,\epsilon ]$ and $\Gamma _2$ is the connected component of $%
B^{-}$ that contains the ``highest'' point of $B^{-}$ in $\{c\}\times
[1-\epsilon ,1].$ So $\Gamma _2$ must also be a non-injective component. $%
\Box $

\vskip 0.2truecm

Summarizing, the above results prove that, for every vertical segment $u$ of
the form $\{l\}\times [0,\epsilon ]$ (or $\{l\}\times [1-\epsilon ,1]$)
which intersects $B^{-}$ (see the definition of $\epsilon >0$ in the
beginning of this section), the ''lowest'' (or ''highest'') component of $%
B^{-}$ in $u$ must be non-injective, $\widetilde{f}$-positively invariant
and dense when projected to $A.$

\section{Proof of theorem 5}

Without loss of generality, suppose $\Gamma \subset B^{-}$ is an injective
down connected component. Consider a vertical $v=\{c\}\times [0,\epsilon [,$
such that:

\begin{equation}
\label{defvene} 
\begin{array}{l}
1) 
\text{ }0<\epsilon \leq \epsilon ^{\prime },\text{ where }\epsilon ^{\prime
} \text{ is such that }\forall (\widetilde{x},\widetilde{y})\in {\rm I}%
\negthinspace {\rm R}\times [0,\epsilon ^{\prime }],\text{ } \\ p_1\circ
\tilde f(\tilde x, 
\widetilde{y})>\tilde x+\sigma ,\text{ for some }\sigma >0;\text{ } \\ 2) 
\text{ }v\subset \Gamma _{down}; \\ 3)\text{ }v\cap B^{-}\neq \emptyset . 
\end{array}
\end{equation}

The above is possible because $\overline{p(B^{-})}=A$ and if $\{c\}\times
[0,b[\subset \Gamma _{down}$ and $(c,b)\in \Gamma ,$ then $B^{-}\cap
\{c\}\times [0,b[\neq \emptyset ,$ see the proof of the previous theorem. So
we can choose a sufficiently small $c$ such that $\{c\}\times [0,\epsilon
^{\prime }]\cap B^{-}\neq \emptyset ,$ where $\epsilon ^{\prime }$ comes
from 1) of (\ref{defvene}). If $\{c\}\times [0,\epsilon ^{\prime }[\subset
\Gamma _{down},$ we are done. If not, let $0<\epsilon <\epsilon ^{\prime }$
be such that $v=\{c\}\times [0,\epsilon [\subset \Gamma _{down}$ and $%
(c,\epsilon )\in \Gamma .$

Denote by $\Theta $ the lowest connected component of $B^{-}$ in $v$ and by $%
w=\{c\}\times [0,\delta [\subset v$ the vertical such that $w\cap
B^{-}=\emptyset $ and $(c,\delta )\in \Theta .$ From theorem 4 we know that $%
\Theta $ satisfies the following conditions:

\begin{equation}
\label{propgaman} 
\begin{array}{l}
i) 
\text{ }\Theta \text{ is non-injective;} \\ ii) 
\text{ }\widetilde{f}(\Theta )\subset \Theta ; \\ iii)\text{ }\overline{%
p(\Theta )}=A; 
\end{array}
\end{equation}

If $\Theta $$\prec \Gamma ,$ then proposition \ref{limitado} implies the
existence of a real number $d$ such that $\Theta $$\cap V_d^{-}\subset
\Gamma _{down}.$ As $dist(\Gamma ,{\rm I}\negthinspace {\rm R}\times
\{1\})>0,$ we get that $dist(\Theta ,{\rm I}\negthinspace 
{\rm R}\times \{1\})>0,$ something that contradicts property iii) of
expression (\ref{propgaman}).

So we can assume that $\Gamma \prec \Theta .$ As in most of the previous
results, let $\Omega $ be the connected component of $(\Theta \cup w)^c$
that contains $]-\infty ,c[\times \{0\}.$ We know that 
\begin{equation}
\label{tttt}closure(\Omega )\supset closure(\Theta _{c,down}) 
\end{equation}
and, as $\Gamma \prec \Theta ,$ $\left[ \Gamma \cap V_c^{-}\right] \subset
\Theta _{c,down}.$ So, using expression (\ref{tttt}) we get that $\Gamma
\subset \Omega $ because $\Gamma$ is connected, 
$\Gamma \cap \Omega \neq \emptyset $ and $\Gamma
\cap \partial \Omega \subset \Gamma \cap (\Theta \cup w)=\emptyset .$ The
rest of our proof will be divided in two steps:

\vskip 0.2truecm

{\bf Step 1: }Here we are going to prove that for all integers $n>0$ and $%
k\geq 0,$ $\widetilde{f}^n(\Gamma )+(k,0)$ is disjoint from $\Theta $ and $
\widetilde{f}^n(\Gamma )\subset \Omega $ $\Rightarrow $ $\widetilde{f}%
^n(\Gamma )+(k,0)\prec \Theta .$

\vskip 0.1truecm

As $\Gamma \prec \Theta ,$ we get for any integer $n>0,$ that either $
\widetilde{f}^n(\Gamma )\subset \Theta $ or $\widetilde{f}^n(\Gamma )\cap
\Theta =\emptyset $ and $\widetilde{f}^n(\Gamma )\prec \Theta ,$ which
implies that $\widetilde{f}^n(\Gamma )\subset \Omega $ because $\widetilde{f}%
^n(\Gamma )\subset B^{-}.$ To begin, suppose $\widetilde{f}(\Gamma )\subset
\Theta .$ This means that $\widetilde{f}^{-1}(\Theta )\supset \Gamma $ and
so, if $\widetilde{f}^{-1}(\Theta )\cap V=\emptyset ,$ then $\widetilde{f}%
^{-1}(\Theta )$ is contained in a connected component of $B^{-},$ that is, $
\widetilde{f}^{-1}(\Theta )=\Gamma ,$ a contradiction because $\Gamma $ is
injective and $\Theta $ is not. So, $\widetilde{f}^{-1}(\Theta )\cap V\neq
\emptyset .$ Let $\Gamma ^{\prime }$ be the connected component of $
\widetilde{f}^{-1}(\Theta )\cap V^{-}$ that contains $\Gamma .$ The fact
that $\widetilde{f}^{-1}(\Theta )$ is connected implies that $\Gamma
^{\prime }$ intersects $V$, is contained in $B^{-}$ and contains $\Gamma .$
So, $\Gamma ^{\prime }=\Gamma $ and this is a contradiction because $\Gamma
\subset \Omega$ and $\Omega\cap V=\emptyset$. So, $\widetilde{f}(\Gamma )\cap \Theta =\emptyset $ $%
\Rightarrow $ $\widetilde{f}(\Gamma )\prec \Theta $ $\Rightarrow $ $
\widetilde{f}(\Gamma )\subset \Omega .$

Now note that, $\widetilde{f}(\Gamma )=\Gamma ^{+},$ because, if this is not
the case, then $\widetilde{f}^{-1}(\Gamma ^{+})\supset \Gamma $ is not a
connected component of $B^{-},$ so $\widetilde{f}^{-1}(\Gamma ^{+})\cap
V\neq \emptyset ,$ which means that $\widetilde{f}^{-1}(\Gamma ^{+})\cap
w\neq \emptyset $ because $\widetilde{f}^{-1}(\Gamma ^{+})\cap \Theta
=\emptyset $ and $\widetilde{f}^{-1}(\Gamma ^{+})\cap \Omega \neq \emptyset
. $ If we denote by $\Gamma ^{+*}$ the connected component of $\widetilde{f}%
^{-1}(\Gamma ^{+})\cap \Omega $ that contains $\Gamma $$,$ then as $
\widetilde{f}^{-1}(\Gamma ^{+})$ is connected, $\Gamma ^{+*}$ intersects $w$
and is contained in $B^{-},$ a contradiction.

So, $\widetilde{f}(\Gamma )=\Gamma ^{+}\subset \Omega $ and an induction
using the above argument implies that for every integer $n>0:$%
$$
\begin{array}{l}
1) 
\text{ }\widetilde{f}^n(\Gamma )\cap \Theta =\emptyset ; \\ 2) 
\text{ }\widetilde{f}^n(\Gamma )\text{ is a connected component of }B^{-};
\\ 3) 
\text{ }\widetilde{f}^n(\Gamma )\subset \Omega ; \\ 4)\text{ }\widetilde{f}%
^n(\Gamma )\prec \Theta ; 
\end{array}
$$

As $\Gamma \subset B^{-}$ is an injective down connected component, the same
holds for $\widetilde{f}^n(\Gamma )$ (for any integer $n>0$). So the
assertion from step 1 holds.

\vskip 0.2truecm

{\bf Step 2: }Here we perform the same construction as we did in lemma \ref
{omegalim}, see it for more details.

\vskip 0.2truecm

Let us fix some $k^{\prime }>0$ in a way that $\widetilde{f}(\Gamma
)+(k^{\prime },0)$ intersects $w.$ Denote by $\Gamma ^{*}$ an unlimited
connected component of $\widetilde{f}(\Gamma )+(k^{\prime },0)\cap
closure(\Omega ).$ By the choice of $k^{\prime }>0$ and the connectivity of $
\widetilde{f}(\Gamma )+(k^{\prime },0),$ we get that $\Gamma ^{*}$ is not
contained in $B^{-}$ because it intersects $v.$ So, there exists a positive
integer $a_1>0$ such that $\widetilde{f}^{a_1}(\Gamma ^{*})$ intersects $%
]0,+\infty [\times [0,1].$

Step 1 implies that $\widetilde{f}^{a_1}(\widetilde{f}(\Gamma )+(k^{\prime
},0))=\widetilde{f}^{a_1+1}(\Gamma )+(k^{\prime },0)$ does not intersect $%
\Theta $ and is smaller then it the order $\prec .$ So, as

$$
\Gamma ^{*}\subset \widetilde{f}(\Gamma )+(k^{\prime },0)\ \Rightarrow \text{
}\widetilde{f}^{a_1}(\Gamma ^{*})\cap \Theta =\emptyset \text{ and } 
\widetilde{f}^{a_1}(\Gamma ^{*})\prec \Theta . 
$$
The fact that $\widetilde{f}^{a_1}(\Gamma ^{*})$ intersects $]0,+\infty
[\times [0,1]$ implies that $\widetilde{f}^{a_1}(\Gamma ^{*})\cap
closure(\Omega )$ has an unlimited connected component, $\Gamma ^{**}$ which
intersects $w.$ So, $\Gamma ^{**}$ is not contained in $B^{-}$ and thus
there exists an integer $a_2>0$ such that $\widetilde{f}^{a_2}(\Gamma
^{**})\subset \widetilde{f}^{a_2+a_1+1}(\Gamma )+(k^{\prime },0)$ intersects 
$]0,+\infty [\times [0,1].$ In exactly the same way as above, we obtain an
unlimited connected component of $\widetilde{f}^{a_2}(\Gamma ^{**})\cap
closure(\Omega )$, denoted $\Gamma ^{***}$ which intersects $w.$ So, $\Gamma
^{***}$ is not contained in $B^{-}$ and there exists an integer $a_3>0$ such
that $\widetilde{f}^{a_3}(\Gamma ^{***})$ intersects $]0,+\infty [\times
[0,1]$ and so on.

Thus, if we define $n_i=a_1+a_2+...+a_i+1,$ we get that $n_i\stackrel{%
i\rightarrow \infty }{\rightarrow }\infty $ and for all $i\geq 1,$ $
\widetilde{f}^{n_i-1}(\widetilde{f}(\Gamma )+(k^{\prime },0))\supset 
\widetilde{f}^{n_i-1}(\Gamma ^{*})=\widetilde{f}^{a_2+...+a_i}(\widetilde{f}%
^{a_1}(\Gamma ^{*}))\supset \widetilde{f}^{a_2+...+a_i}(\Gamma ^{**})\supset
...\supset \widetilde{f}^{a_i}(\Gamma ^{\stackrel{i-times}{*...*}})$ and $
\widetilde{f}^{a_i}(\Gamma ^{\stackrel{i-times}{*...*}})$ intersects $%
]0,+\infty [\times [0,1].$ So, 
$$
\widetilde{f}^{n_i}(\Gamma )\text{ intersects }V_0-(k^{\prime
},0)=V_{-k^{\prime }} 
$$
and this contradicts theorem 2 and thus proves theorem 5.

\vskip0.2truecm

{\it Acknowledgements: }The second author would like to thank Professor
Patrice Le Calvez for telling him to study the set $B^{-}.$

\centerline{\bf Figure captions.}

\begin{itemize}
\item[Figure 1. ]  Diagram showing $\widehat{A}.$

\item[Figure 2. ]  Diagram showing the set $\Gamma_N.$

\item[Figure 3. ]  Diagram showing that 
$\left[ \Gamma _2\cap V_a^{-}\right] \subset 
\Gamma_{1a,down}\cup \Gamma _{1a,up}.$

\item[Figure 4. ]  Diagram showing that either 
$\left[ \Gamma _2\cap V_a^{-}\right] \subset \Gamma _{1a,down}$ or $\left[
\Gamma _2\cap V_a^{-}\right] \subset \Gamma _{1a,up}.$

\item[Figure 5. ]  Diagram showing the sets $\Gamma$, $\Gamma-(1,0)$ and 
$\gamma.$

\item[Figure 6. ]  Diagram showing the sets 
$\overline{p(\Gamma )}\subset \gamma _E^{-}.$ 

\end{itemize}

\begin{center}
\mbox{\includegraphics[width=11cm]{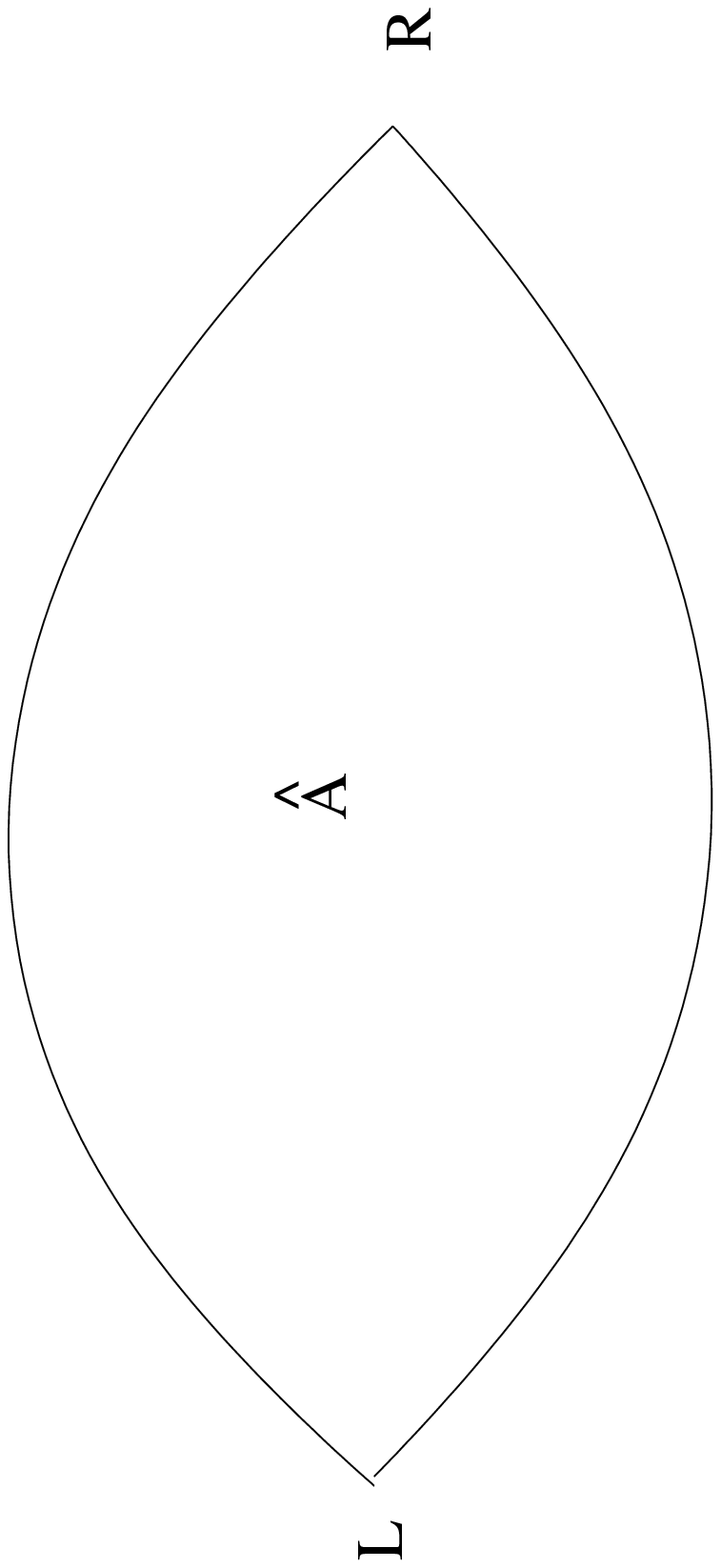}}

\end{center}

\begin{center}
\mbox{\includegraphics[width=11cm]{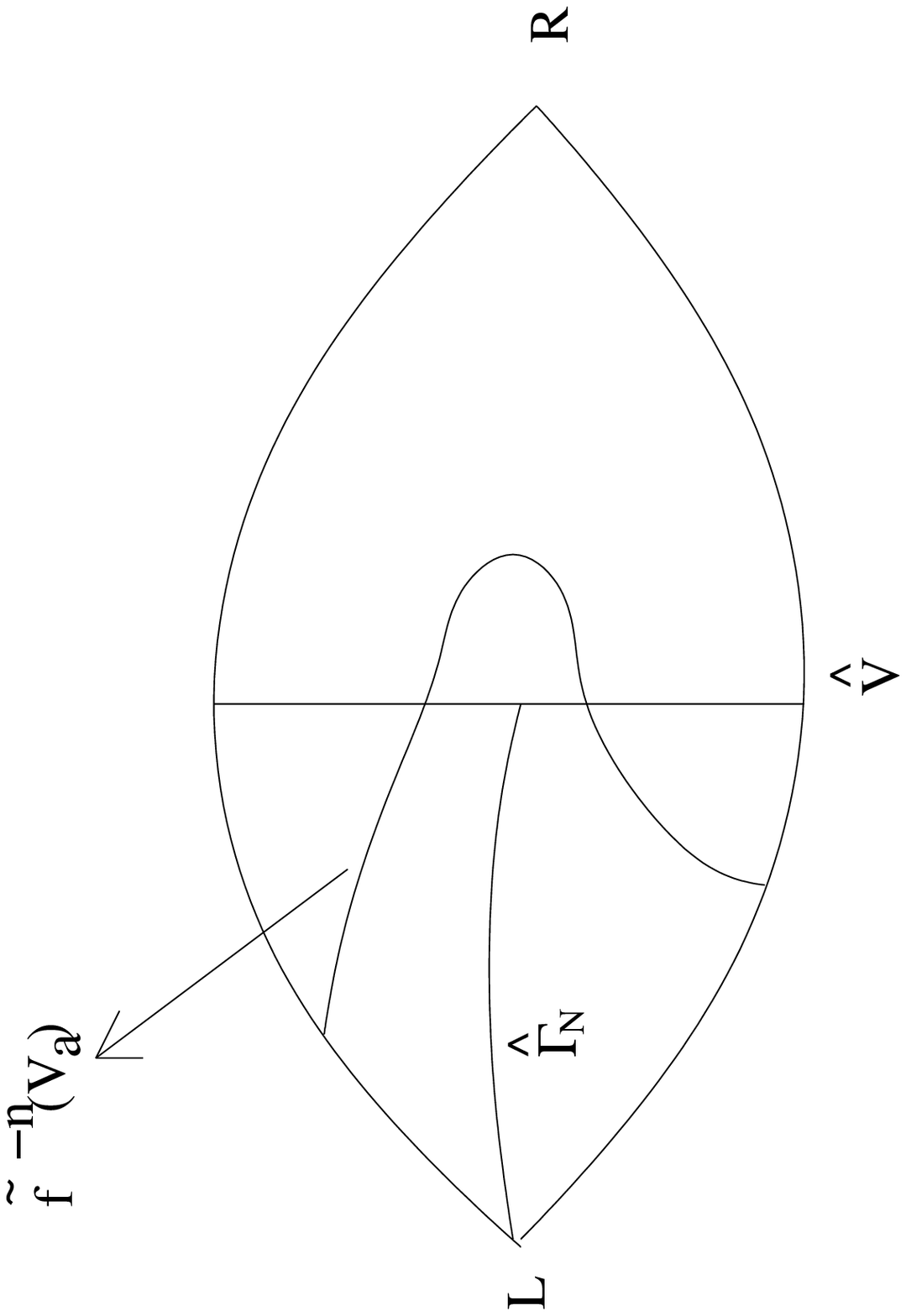}}

\end{center}

\begin{center}
\mbox{\includegraphics[width=11cm]{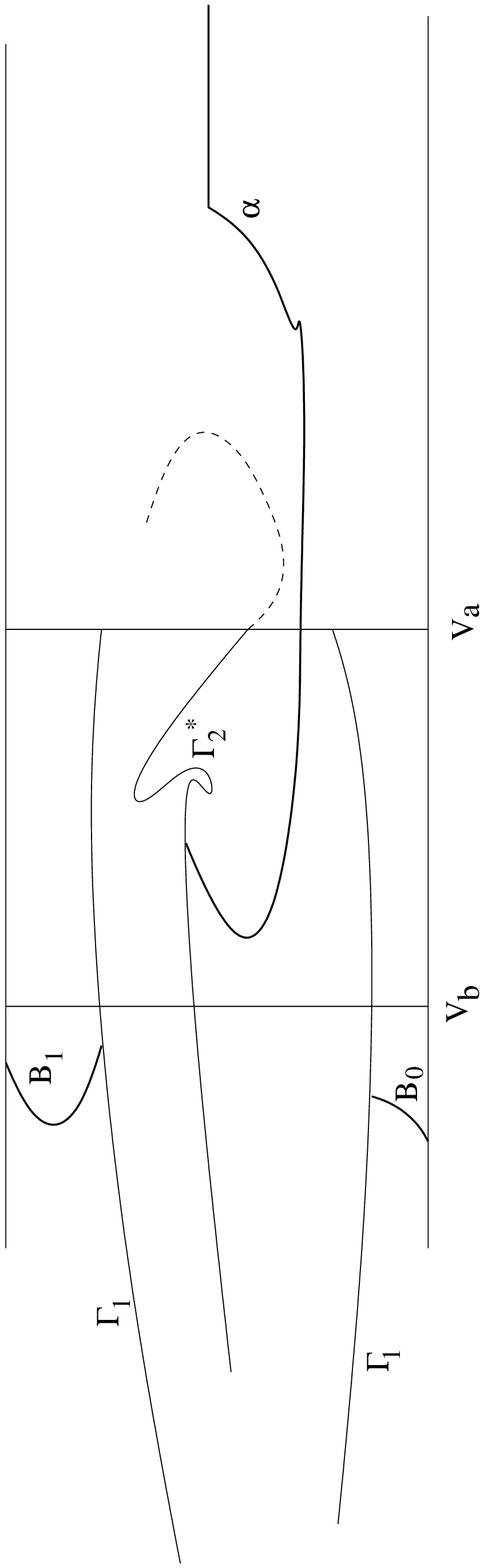}}

\end{center}

\begin{center}
\mbox{\includegraphics[width=11cm]{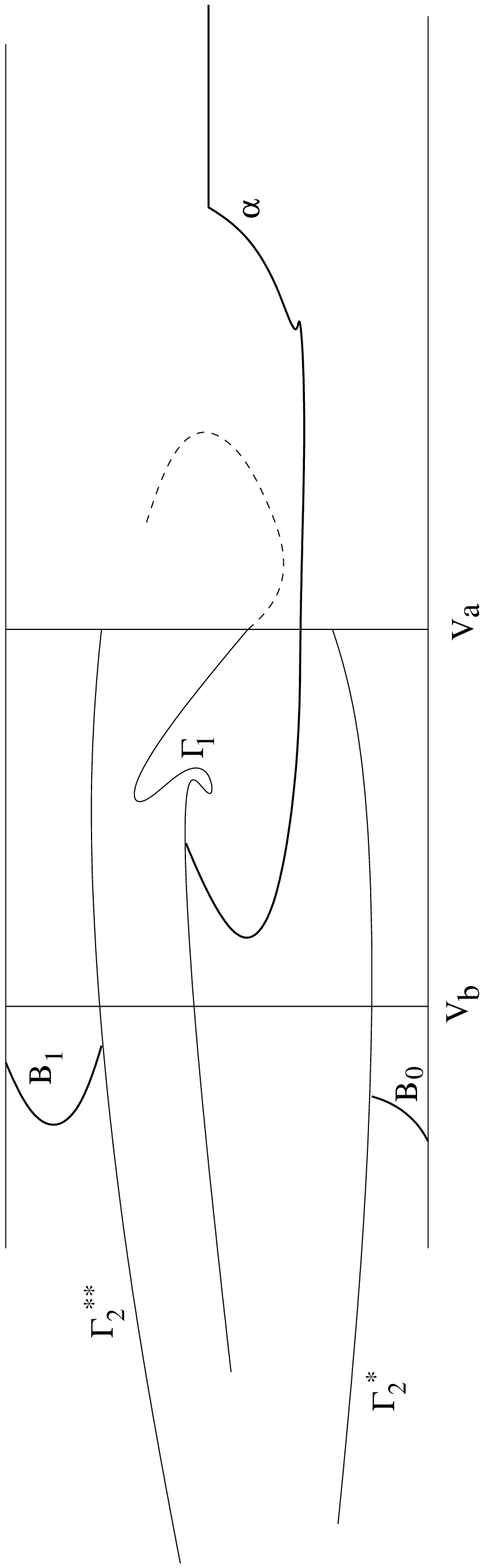}}

\end{center}

\begin{center}
\mbox{\includegraphics[width=11cm]{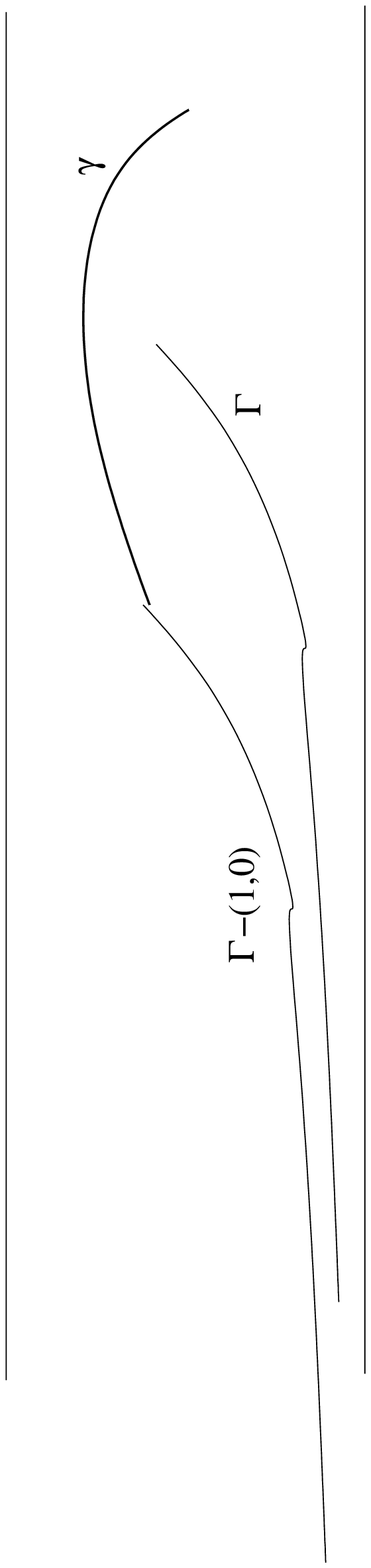}}

\end{center}

\begin{center}
\mbox{\includegraphics[width=11cm]{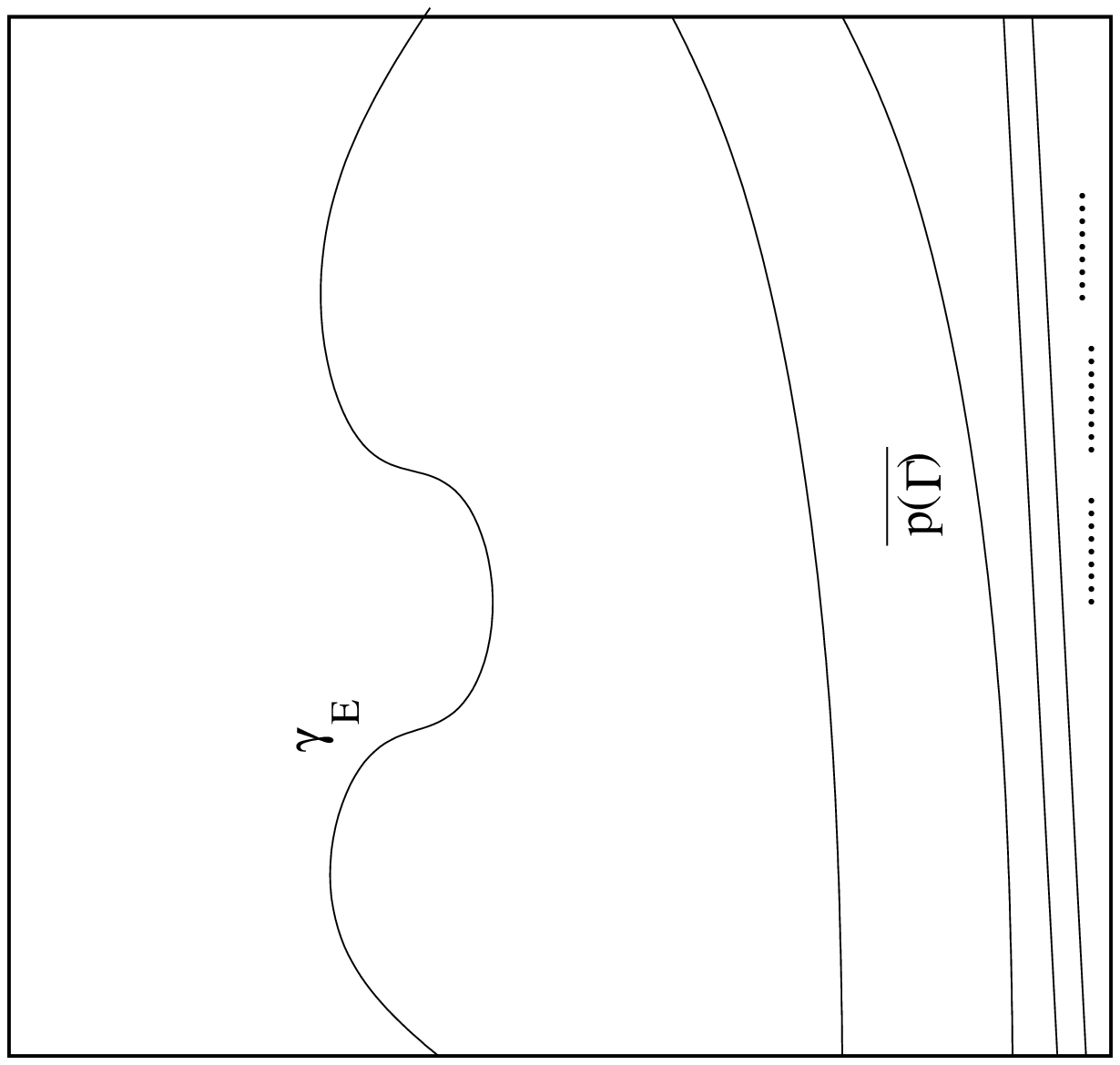}}

\end{center}

\end{document}